\documentclass[a4paper,12pt,leqno,twoside]{article}
\usepackage[english]{babel} 
\usepackage{amsmath}
\usepackage{amsfonts}
\usepackage{amssymb}
\usepackage{latexsym}
\usepackage{xypic}
\usepackage[all]{xy}
\usepackage{stmaryrd} 
\usepackage{theorem}  
\usepackage[mathscr]{eucal}
\usepackage{eufrak}
\usepackage{epsf}
\usepackage[dvips]{graphicx}

\usepackage{epsfig}

\usepackage{latexsym}  
\usepackage[usenames]{color}

\usepackage{fancyhdr}
\usepackage{mathrsfs} 

\usepackage{helvet}

\newcommand{\mc}[1]{\mathcal{#1}}
\newcommand{\ms}[1]{\mathscr{#1}}
\newcommand{\msf}[1]{\mathsf{#1}}
\renewcommand{\phi}{\varphi}
\renewcommand{\theta}{\vartheta}
\renewcommand{\rho}{\varrho}
\newcommand{\ep}{\epsilon}

\newcommand{\lra}{\longrightarrow}

\renewcommand{\bigr}[1]{{\big(#1\big)}}
\renewcommand{\Bigr}[1]{{\Big(#1\Big)}}
\renewcommand{\biggr}[1]{{\bigg(#1\bigg)}}
\newcommand{\bigc}[1]{{\big\{#1\big\}}}
\newcommand{\Bigc}[1]{{\Big\{#1\Big\}}}

\newcommand{\ou}[3][]{\overset{{#1}}{\underset{{#2}}{{#3}}}}

\newcommand{\com}{\mathbb{C}}
\newcommand{\rea}{\mathbb{R}}
\newcommand{\rat}{\mathbb{Q}}
\newcommand{\integer}{\mathbb{Z}}
\newcommand{\integergz}{\mathbb{Z}_{>0}}

\newcommand{\reagz}{\mathbb{R}_{>0}}

\newcommand{\re}{\mathrm{Re}\,}

\newcommand{\zcz}{z^{-1}\com[z^{-1}]}

\newcommand{\zlczl}{z^{-1/l}\com[z^{-1/l}]}

\newcommand{\dist}[2]{\mathrm{dist}(#1,#2)}

\newcommand{\csect}[3]{{S_{#1\pm #2,#3}}}

\newcommand{\Op}{\mathrm{Op}}

\newcommand{\opxsac}{\mathrm{Op}^c(\xsa)}

\newcommand{\xsa}{{X_{sa}}}

\newcommand{\yxsa}{{Y_{\xsa}}}

\newcommand{\Ho}[3][]{\mathcal{H}\mathrm{om}_{#1}(#2,#3)}

\newcommand{\hob}[3][]{\mathrm{Hom}_{#1}\big(#2,#3\big)}

\newcommand{\M}{\mathcal{M}}
\newcommand{\ot}{\mathcal{O}^t}
\newcommand{\otxsa}{\mathcal{O}^t_\xsa}

\newcommand{\dbt}{\mathcal{D}b^t}

\newcommand{\dbxr}{\mathcal{D}b_{X_\rea}}

\newcommand{\dbtxsa}{\mathcal{D}b^t_{\xsa}}

\newcommand{\D}{\mathcal{D}}

\renewcommand{\mod}{\mathrm{Mod}}
\renewcommand{\M}{\mathcal{M}}

\newcommand{\st}[1]{\mathscr{S}^t\big(#1\big)}

\newcommand{\sh}[1]{\mathscr{S}\big(#1\big)}

\newtheorem{thm}{Theorem}[subsection]
\newtheorem{df}[thm]{Def\mbox{}inition}
\newtheorem{cor}[thm]{Corollary}
\newtheorem{prop}[thm]{Proposition}
\newtheorem{lem}[thm]{Lemma}

\newtheorem{rem}[thm]{Remark}

\numberwithin{equation}{section}

\newcommand{\ben}{\begin{enumerate}}
\newcommand{\een}{\end{enumerate}}

\newcommand{\proofend}{\hfill $\Box$ \vspace{\baselineskip}\newline}
\author{Giovanni Morando}
\title{\bf{Tempered solutions of $\D$-modules
    on complex curves and formal invariants}}
\date{December 5, 2007}

\long\def\symbolfootnote[#1]#2{\begingroup%
\def\thefootnote{\fnsymbol{footnote}}\footnote[#1]{#2}\endgroup}

\newpage

\begin{document}

\maketitle

\thispagestyle{empty}

\begin{abstract}
\sloppy 
Let $X$ be a complex analytic curve. In this paper we prove that the
subanalytic sheaf of tempered holomorphic solutions of $\D_X$-modules
induces a fully faithful functor on a subcategory of germs of formal
holonomic $\D_X$-modules. Further, given a germ $\mc M$ of holonomic
$\D_X$-module, we obtain some results linking the subanalytic sheaf of
tempered solutions of $\mc M$ and the classical formal and analytic
invariants of $\mc M$.  
\end{abstract}

\symbolfootnote[0]{\phantom{a}\hspace{-7mm}\textit{2000 MSC.} Primary
  34M35; Secondary 32B20 34Mxx.} 

\symbolfootnote[0]{\phantom{a}\hspace{-7mm}\textit{Keywords and
    phrases:} $\D$-modules, irregular singularities, tempered
  holomorphic functions, subanalytic.} 

\vspace{-13mm}

\tableofcontents

\phantom{a}

\phantom{a}


\section*{Introduction}\markboth{Introduction}{Introduction}
\addcontentsline{toc}{section}{\textbf{Introduction}}

The search for algebraic or topological invariants of complex linear
partial dif\mbox{}ferential equations is classical and widely developed. 

At the very f\mbox{}irst step of the study of linear dif\mbox{}ferential equations,
two main types of equations are distinguished: regular and
irregular. To give an idea of the dif\mbox{}ference between the two kinds of
equations, let us recall that, in dimension $1$, the solutions
of the former equations have moderate growth while the solutions of
the latter have exponential-type growth. 

The more general algebraic approach to the study of linear
dif\mbox{}ferential equations consists in considering dif\mbox{}ferential equations as sheaves of
modules over the ring $\D_X$ of linear dif\mbox{}ferential operators on a
manifold $X$. In this 
framework, in \cite{kashi79} and \cite{rims}, M. Kashiwara gives a proof of the
Riemann-Hilbert correspondence which is a generalization of the 21st
Hilbert's problem. For $X$ a complex analytic manifold, M. Kashiwara def\mbox{}ines the
functor $T\mc Hom$ and he gives an explicit inverse to the functor of
holomorphic solutions from 
the bounded derived category of complexes of $\D_X$-modules with
regular holonomic cohomology to 
the bounded derived
category of complexes of sheaves with constructible cohomology. This
implies the classic result that the functor of holomorphic solutions
$\ms S(\cdot)$ is an equivalence between the category of regular
connections on $X$ with poles 
along a closed submanifold $Z$ and the category of linear
representations of f\mbox{}inite dimension of the fundamental group of
$X\setminus Z$. 

The irregular case is more complicated. In complex dimension 1, the classif\mbox{}ication of
meromorphic connections through the formal classif\mbox{}ication and the
Stokes phenomena is nowadays well understood. 
 Let us roughly explain the formal
classif\mbox{}ication. It is based on the Levelt-Turrittin's Formal
Theorem. Such theorem is one of the cornerstones of the study of
ordinary dif\mbox{}ferential equations. For $\phi\in\zcz$, set $\mc
L^\phi:=\D_\com\exp(\phi)$. The $\D_\com$-module $\mc L^\phi$ is
one of the basic examples   
 of an ir\-re\-gu\-lar meromorphic connection. 
Given a f\mbox{}inite set $\Sigma\subset\zcz$, a family of regular
meromorphic connections indexed by $\Sigma$, $\{\mc
R_\phi\}_{\phi\in\Sigma}$, one calls $\ou{\phi\in\Sigma}{\oplus}{\mc
L^{\phi}\otimes\mc R_\phi}$ a good model. Let $\mu_l:\com\to\com$,
$z\mapsto z^l$ and $\mu_l^*$ the pull-back functor for germs of $\D$-modules. The
Levelt-Turrittin's Formal Theorem states that, given a meromorphic
connection $\mc M$, there exists $l\in\integergz$ such that
\begin{equation}\label{eq_intro_fht}
\mu_l^*\mc M\simeq \ou{\phi\in\Sigma}{\oplus}{\mc 
L^{\phi}\otimes\mc R_\phi} \ ,
\end{equation}
as modules over the ring of linear
dif\mbox{}ferential operators with formal Laurent series
coef\mbox{}f\mbox{}icients.  
The formal isomorphism given in \eqref{eq_intro_fht} does not induce an isomorphism in the analytic
category, that is to say we can have two non-isomorphic meromorphic
connections which are formally isomorphic. Such peculiarity is at the
base of the Stokes phenomena. Roughly speaking, the
Hukuhara-Turrittin's Asymptotic Theorem states that the isomorphism
\eqref{eq_intro_fht} is analytic on suf\mbox{}f\mbox{}iciently small
open sectors. Indeed, such a theorem states that the holomorphic
solutions of a meromorphic connection $\mc M$ on a small open sector $S$ are
$\com$-linear combinations of a f\mbox{}inite number of functions of
the form $h(z)\exp(\phi(z^{1/l}))$, for $z^{1/l}$ an arbitrary $l$-th
root of $z$ defined on $S$, $\phi\in\zcz$ not depending on $S$ and $h$ an invertible
tempered holomorphic function on $S$ with tempered inverse. The
functions $\phi(z^{1/l})$, which do not depend on $S$, are called
determinant polynomials of $\mc M$.

In this paper we will use the description of meromorphic connections
given by $\Omega$-f\mbox{}iltrations, introduced by P. Deligne in
\cite{dmr}. Let us brief\mbox{}ly recall the setting of
$\Omega$-f\mbox{}iltrations. We refer to \cite{dmr}, \cite{malgr_class}, \cite{bv} and \cite{malgr_birk}. First one def\mbox{}ines
a local system $\Omega$ on $S^1$ whose stalks are partially ordered
sets. Such local system represents the space of all possible
determinant polynomials. Further, one def\mbox{}ines the category of
$\Omega$-f\mbox{}iltered and $\Omega$-graded local systems. Roughly speaking,
$\Omega$-f\mbox{}iltered local systems are obtained by gluing $\Omega$-graded
local systems preserving the local order of $\Omega$. Now, the
holomorphic solutions of a meromorphic connection can be endowed with
a structure of $\Omega$-f\mbox{}iltered local system. Hence one gets a
functor, $\ms S^\Omega$, which is an equivalence of categories between
meromorphic connections (resp. formal meromorphic connections) and
$\Omega$-f\mbox{}iltered (resp. $\Omega$-graded) local systems. Roughly speaking, given a
meromorphic connection $\mc M$ and $\theta\in S^1$, the $\Omega$-graduation on $\ms
S^\Omega(\mc M)_\theta$ is determined by the positive integer $l$, the
set $\Sigma$ and the rank of the local systems $\ms S(\mc R_\phi)$ appearing in \eqref{eq_intro_fht}. The gluing
morphisms def\mbox{}ining the $\Omega$-f\mbox{}iltration on $\ms 
S^\Omega(\mc M)$ are determined by the Stokes phenomena and the formal
monodromy (i.e. the monodromy of the local systems $\ms S(\mc
R_\phi)$). Even though the $\Omega$-f\mbox{}iltrations provide a complete description of the category of
meromorphic connections, the local system $\Omega$ remains an analytic object
built ad hoc. It is interesting to look for a topological description
of the determinant polynomials. 

Let us also recall that many sheaves of function spaces have been used in the
study of irregular ordinary dif\mbox{}ferential equations. For example, one
can f\mbox{}ind in \cite{malgr_birk} the def\mbox{}initions of the sheaves $\mathcal
A^{\leq r}$ ($r\in\rea$) def\mbox{}ined on the real blow-up of the complex plane at the
origin. In \cite{dmr}, P. Deligne def\mbox{}ined the sheaves
$\widetilde{\mathcal F}^{k}$, successively studied in detail in \cite{loday_pourcin}. Roughly
speaking, the solutions of $\mc L^\phi$ with values in $\mc A^{\leq
r}$ (resp. $\widetilde{\mathcal F}^{k}$) depend only on the degree and
the argument of the leading coef\mbox{}f\mbox{}icient of $\phi$ (resp. the degree and the
leading coef\mbox{}f\mbox{}icient of $\phi$).

In higher dimension the study of irregular $\D$-modules is much more
complicated. In \cite{sabbah_ast}, C. Sabbah def\mbox{}ines the notion of
good model in dimension $2$, he conjectures the analogue
of the formal Levelt-Turrittin's Theorem and he proves it for
meromorphic connections of rank $\leq5$. Recently, T. Mochizuki
announced the proof of Sabbah's conjecture for the algebraic case. 


Given a complex analytic manifold $X$, in \cite{ast}, M. Kashiwara and
P. Schapira def\mbox{}ined the complex of sheaves of tempered holomorphic
functions $\otxsa$. The entries of the complex $\otxsa$ are not
sheaves on a topological space but on the subanalytic site,
$\xsa$. The open sets of $\xsa$ are the subanalytic open subsets of
$X$, the coverings are locally f\mbox{}inite coverings. If $X$ has dimension $1$, then $\otxsa$ is
a sheaf on $\xsa$ and, for $U$ a relatively compact subanalytic open
subset of $X$, the sections of $\otxsa(U)$ are the holomorphic
functions on $U$ which extend as distributions on $X$ or,
equivalently, which have moderate growth at the boundary of $U$. 

Further in an example in \cite{arx}, M. Kashiwara and P. Schapira explicited 
the sheaf of tempered holomorphic solutions of $\mc L^{1/z}$. Such example suggests that tempered
holomorphic functions and the subanalytic site could be useful tools
in the study of ordinary dif\mbox{}ferential equations. Roughly speaking,
Deligne's idea of $\Omega$-f\mbox{}iltrations aimed to enrich the structure
of the category of sheaves of $\com$-vector spaces where the functor of
holomorphic solutions takes values. The approach
through subanalytic sheaves enriches the topology of the space where
the sheaf of solutions lives. 

In this paper we go into the study of the subanalytic sheaf of
tempered holomorphic solutions of germs of $\D$-modules. Denote by
$\st {\mc M}$ the subanalytic sheaf of tempered 
holomorphic solutions of a holonomic $\D_X$-module $\mc M$. Let
$X\subset\com$ be an open neighborhood of $0$, $\mod(\com_\xsa)$ the
category of sheaves of $\com$-vector spaces on $\xsa$, $\mathsf
{GM}_k$ be the subcategory 
of germs at $0$ of $\D$-modules consisting of good models with Katz
invariant $<k$. We prove that 
$$ \st{\cdot\otimes\mc L^{1/z^k}}:\mathsf{GM}_k\lra\mod(\com_\xsa) $$
is a fully faithful functor ({\bf Theorem \ref{thm_ff}}). Further we
prove that, given a germ of holonomic $\D$-module $\mc M$ with Katz
invariant $<k$, the datum of $ \st{\mc M\otimes\mc L^{1/z^{k+1}}}$ 
is equivalent to the data of the stalk of $\ms S^\Omega(\mc M)$ and
the local system underlying $\ms S^\Omega(\mc M)$ (i.e. the local
system of holomorphic solutions) ({\bf Theorem \ref{thm_final}}).  

In conclusion we can say that tempered solutions on the subanalytic
site give a topological description of the determinant polynomials
of a given meromorphic connection. As further developement, it would
be interesting to describe precisely the image category of the functor of
tempered solutions in order to give a topological description of the
space of determinant polynomials. It would also be interesting to give
a good notion of Fourier trasform for tempered holomorphic solutions
of algebraic $\D$-modules in the same spirit of \cite{malgr_birk}. 

The present paper is subdivided in three sections organized as
follows. 

{\bf Section 1} is devoted to the def\mbox{}initions, the notations and the
presentation of the main results that will be needed in the rest of
the paper. In particular we recall classical results on the subanalytic
sets and site, the tempered holomorphic functions and the germs of
$\D$-modules.

The functions of the form $\exp(\phi)$, $\phi\in\zcz$, are the
responsible for the non-tempered-growth of the solutions of an
irregular $\mc D$-module. This motivates the study of $\exp(\phi)$
that we develop in {\bf Section 2}. In particular, given
$\phi\in\zcz$ and $U$ a relatively compact subanalytic open subset of $\com$, we give a
necessary and suf\mbox{}f\mbox{}icient topological condition on $U$ so that
$\exp(\phi)\in\ot_{\com_{sa}}(U)$. Further, we prove that the
condition ``for any $U\subset\com$ relatively compact subanalytic
open set, $\exp(\phi_1)\in\ot_{\com_{sa}}(U)$ if and only if
$\exp(\phi_2)\in\ot_{\com_{sa}}(U)$'' is 
equivalent to ``$\phi_1$ and $\phi_2$ are proportional by a real
positive constant''. 

In {\bf Section 3} we apply the results of
Section 2 to the study of the functor of tempered holomorphic solutions of
germs of $\D$-modules on a complex curve. We prove that $
\st{\cdot\otimes\mc L^{1/z^k}}:\mathsf{GM}_k\to\mod(\com_\xsa) $ is a
fully faithful functor and that, given a germ of $\D$-module $\mc M$ with Katz
invariant $<k$, the datum of $ \st{\mc M\otimes\mc L^{1/z^{k+1}}}$ is
equivalent to the data of the determinant polynomials of $\mc M$ and
holomorphic solutions of $\mc M$.  

\emph{Aknowledgements:} We thank P. Schapira for proposing this
problem to our attention, C. Sabbah and N. Honda for many fruitful
discussions and A. D'Agnolo for many useful remarks.

\section{Notations and review}  \label{SEC_AN_GEOM}

In this section we recall the def\mbox{}initions and the
classical results concerning:
\begin{enumerate}
\item subanalytic sets, the subanalytic site
and sheaves on it,
\item  the subanalytic sheaf of tempered holomorphic
functions,
\item germs of $\D$-modules on a complex curve.
\end{enumerate}

\subsection{The subanalytic site}\label{SUBSEC_RECALL_XSA}

Let $M$ be a real analytic manifold, $\mathcal A$ the sheaf of
real-valued real analytic functions on $M$.

\begin{df}
\begin{enumerate}
\item A set $X\subset M$ is said \emph{semi-analytic at} $x\in M$ if the following condition is satisf\mbox{}ied. There exists an open neighborhood $W$ of $x$ such that $X\cap
 W=\cup_{i\in I}\cap_{j\in J}X_{ij}$ where $I$ and $J$ are f\mbox{}inite sets
 and either $X_{ij}=\{y\in W;\ f_{ij}(y)>0\}$ or $X_{ij}=\{y\in W;\
 f_{ij}(y)=0\}$ for some $f_{ij}\in\mathcal A(W)$. Further, $X$ is said \emph{semi-analytic} if $X$ is semi-analytic at any $x\in M$. 
\item A set $X\subset M$ is said \emph{subanalytic} if the following condition is satisf\mbox{}ied. For any $x\in M$, there exist an open neighborhood $W$ of $x$, a real analytic manifold $N$ and a
  relatively compact semi-analytic set $A\subset M\times N$ such that
  $\pi(A)=X\cap W$, where $\pi: M\times N\to M$ is the
  projection. 
\end{enumerate}
\end{df}

Given $X\subset M$, denote by $\overset{\circ}{X}$ (resp. $\overline
X$, $\partial X$) the interior (resp. the closure, the boundary) of $X$. 

\begin{prop}[See \cite{bier_mil}]
Let $X$ and $Y$ be subanalytic subset of $M$. Then $X\cup Y$, $X\cap
Y$, $\overline X$, $\overset{\circ}{X}$ and $X\setminus Y$ are
subanalytic. Moreover the connected components of $X$ are subanalytic,
the family of connected components of $X$ is locally f\mbox{}inite and $X$ is
locally connected.
\end{prop}

Proposition \ref{prop_reg_dist} below is based on \L ojasiewicz's
inequality, see \cite[Corollary 6.7]{bier_mil}. 

\begin{prop}\label{prop_reg_dist}
Let $U\subset\rea^n$ be an open set, $X,Y$ closed subanalytic subsets
of $U$. For any $x_0\in X\cap Y$, there exist an open neighborhood $W$
of $x_0$, $c,r\in\reagz$ such that, for any $x\in W$,
$$ \dist xX +\dist xY\geq c\,\dist x{X\cap Y}^r\ .$$
\end{prop}


\begin{df}
Let $\ep\in\reagz$, $\gamma:]-\ep,\ep[\to M$ an analytic map. The set
$\gamma(]0,\ep[)$ is said a \emph{semi-analytic arc with an endpoint at} $\gamma(0)$.
\end{df}

\begin{thm}[Curve Selection Lemma.]\label{thm_csl}
Let $Z\neq\varnothing$ be a subanalytic subset of $M$ and let
$z_0\in\overline Z$. Then there exists an ana\-ly\-tic map $
\gamma:]-1,1[\longrightarrow M $, such that $\gamma(0)=z_0$ and
$\gamma(t)\in Z$ for $t\neq0$. 
\end{thm}

For the rest of the subsection we refer to \cite{ast}.

Let $X$ be a complex analytic curve, we denote by $\Op(X)$ the family
of open subsets of $X$. For $k$ a commutative unital ring, we
denote by $\mod(k_X)$ the category of sheaves of $k$-modules on $X$.

Let us recall the def\mbox{}inition of the subanalytic site $\xsa$ associated
to $X$. An element $U\in\Op(X)$ is an open set for $\xsa$ if it is open,
relatively compact and subanalytic in $X$. The family of open sets of $\xsa$
is denoted $\opxsac$. For $U\in\opxsac$, a subset $S$ of the family of
open subsets of $U$ is said an open covering of $U$ in $\xsa$ if
$S\subset\opxsac$ and, for any compact $K$ of
$X$, there exists a f\mbox{}inite subset $S_0\subset S$ such that
$K\cap(\cup_{V\in S_0}V)=K\cap U$. 

We denote by $\mod(k_\xsa)$ the category of sheaves of $k$-modules on the
subanalytic site. With the aim of def\mbox{}ining the category $\mod(k_\xsa)$, the adjective ``relatively compact'' can be omitted in the def\mbox{}inition above. Indeed, in \cite[Remark 6.3.6]{ast}, it is proved that
$\mod(k_\xsa)$ is equivalent to the category of sheaves on the site
whose open sets are the open subanalytic subsets of $X$ and whose
coverings are the same as $\xsa$.

Given $Y\in\opxsac$, we denote by $\yxsa$ the site induced by $\xsa$
on $Y$, def\mbox{}ined as follows. The open sets of $\yxsa$ are open subanalytic subsets
of $Y$. A covering of $U\in\Op(Y_{sa})$ for the topology $\yxsa$ is a
covering of $U$ is $\xsa$. 

We denote by $ \rho:X\lra\xsa $, the natural morphism of sites
associated to $\opxsac\lra\Op(X)$. We refer to \cite{ast} for the
def\mbox{}initions of the functors $\rho_*:\mod(k_X)\lra\mod(k_\xsa)$
and $\rho^{-1}:\mod(k_\xsa)\lra\mod(k_X)$ and for Proposition
\ref{prop_functors} below. 
\begin{prop}\label{prop_functors}
\begin{enumerate}
\item $\rho^{-1}$ is left adjoint to $\rho_*$.
\item $\rho^{-1}$ has a left adjoint denoted by $\rho_!:\mod(k_X)\lra\mod(k_\xsa)$ .
\item $\rho^{-1}$ and $\rho_!$ are exact and $\rho_*$ is exact on $\rea$-constructible sheaves.
\item $\rho_*$ and $\rho_!$ are fully faithful.
\end{enumerate}
\end{prop}

Through $\rho_*$, we will consider $\mod(k_X)$ as a subcategory of $\mod(k_\xsa)$.

The functor $\rho_!$ is described as follows. If
$U\in\opxsac$ and $F\in\mod(k_X)$, then $\rho_!(F)$ is the sheaf on $\xsa$
associated to the presheaf $U\mapsto F\bigr{\overline U}$. 

\begin{rem}
  It is worth to mention that, given an analytic manifold $X$, there
  exists a topological space $X'$ such that the category of sheaves on
  $\xsa$ with values in sets is equivalent to the category of sheaves
  on $X'$ with values in sets. A detailed description of the
  semi-algebraic case and the o-minimal case are presented
  respectively in \cite{bcr} and \cite{coste}. 
\end{rem}

\subsection{Def\mbox{}inition and main properties of
  $\ot_\xsa$}\label{SUBSEC_OT}

For this subsection we refer to \cite{ast}.

Let $X$ be a complex analytic curve, denote by $\overline X$ the
complex conjugate curve and by $X_\rea$ the underlying real analytic
manifold. We denote by $\xsa$ the subanalytic site relative to $X_\rea$.


Denote by $\mathcal O_X$ (resp. $\D_X$) the sheaf of holomorphic
functions (resp. linear dif\mbox{}ferential operators with holomorphic
coef\mbox{}f\mbox{}icients) on $X$. Denote by $\dbxr$ the sheaf of distributions on $X_\rea$ and, for a
closed subset $Z$ of $X$, by $\Gamma_Z(\dbxr)$ the subsheaf of
sections supported by $Z$. One denotes by $\dbtxsa$ the presheaf of
tempered distributions on $\xsa$ def\mbox{}ined as follows,
$$\opxsac\ni U \longmapsto \dbtxsa(U):=\frac{\Gamma(X;\dbxr)}{\Gamma_{X\setminus U}(X;\dbxr)} \ .$$
In \cite{ast}, using some results of \cite{loj_studia}, it is proved that $\dbtxsa$ is a
sheaf on $\xsa$. This sheaf is well def\mbox{}ined in the category
$\mod(\rho_! \D_X)$. 
Moreover, for any $U\in\opxsac$, $\dbtxsa$ is $\Gamma(U,\cdot)$-acyclic.

Denote by $D^b\bigr{\rho_!\D_X}$ the bounded derived category of
$\rho_!\D_X$-modules. The sheaf $\ot_\xsa \in D^b\bigr{\rho_!\D_X}$ 
of tempered holomorphic functions is def\mbox{}ined as
$$\ot_\xsa:=R\mathcal Hom_{\rho_! \mathcal{D}_{\overline X}}\bigr{\rho_!\mathcal{O}_{\overline X},\dbt_{X_{\rea}}}\ .$$

In \cite{ast}, it is proved that, since $\dim X=1$, $R\rho_*\mathcal O_X$
and $\ot_\xsa$ are concentrated in degree $0$ .
Hence we can write the following exact sequence of sheaves on $\xsa$

$$ 0\longrightarrow \ot_\xsa \longrightarrow \dbtxsa
\overset{\bar\partial}{\longrightarrow} \dbtxsa\longrightarrow 0\ .$$

Let us recall that the def\mbox{}inition of $\dbtxsa$ and $\otxsa$
can be given without any change in the case  of $X$ a complex analytic
manifold (see \cite{ast}).

Now we recall the def\mbox{}inition of polynomial growth for
$\mathcal{C}^\infty$ functions on $X_\rea$ and in
\eqref{eq_debar_cinftyt} we give an alternative expression for
$\ot_\xsa(U)$, $U\in\opxsac$. 

\begin{df}
Let $U$ be an open subset of $X$,
$f\in\mathcal{C}^\infty_{X_\rea}(U)$. One says that $f$ has
\emph{polynomial growth at $p\in X$} if it satisf\mbox{}ies the
following condition. For a local coordinate system $x=(x_1,x_2)$
around $p$, there exist a compact neighborhood $K$ of $p$ and $M\in\integergz$ such that 
\begin{equation}\label{eq_temp}
 \underset{x\in K\cap U}{\sup}\dist x{K\setminus U}^M\big|f(x)\big|<+\infty \ .
\end{equation}
We say that $f\in\mathcal{C}^\infty_{X_\rea}(U)$ has {\em polynomial
  growth on $U$} if it has
polynomial growth at any $p\in X$. We say that $f$ is {\em tempered at
  $p$} if all its derivatives have polynomial growth at $p\in X$. We say
that $f$ is {\em tempered on $U$} if it is tempered at any $p\in X$. Denote
  by $\mathcal C^{\infty,t}_{\xsa}$ the presheaf on $X_\rea$ of tempered
  $\mathcal C^{\infty}$-functions.
\end{df}

It is obvious that $f$ has polynomial growth at any point of $U$. 

In \cite{ast} it is proved that $\mathcal C^{\infty,t}_{\xsa}$ is a sheaf
on $\xsa$. For $U\subset\rea^2$ a relatively compact open set, there
is a simple characterization of functions with polynomial growth on
$U$. 

\begin{prop}\label{prop_equiv_def_pol_gr}
Let $U\subset\rea^2$ be a relatively compact open set and let $f\in\mathcal{C}^\infty_{\rea^2}(U)$.
Then $f$ has polynomial growth if and only if  there exist $C,M\in\reagz$ such that, for any $x\in U$, 
$$\big|f(x)\big|\leq\frac C{\dist{x}{\partial U}^{M}} \ .$$
\end{prop}

For Proposition \ref{otct} below, see \cite{ast}. 
\begin{prop}\label{otct}
One has the following isomorphism
$$ \ot_\xsa\simeq R\mathcal Hom_{\rho_!\D_{\overline
    X}}\bigr{\rho_!\mathcal{O}_{\overline
    X},\mathcal{C}^{\infty,t}_{\xsa}}\ .$$ 
\end{prop}

Hence, for $U\in\opxsac$, we deduce the short exact sequence
\begin{equation}\label{eq_debar_cinftyt} 0\longrightarrow \ot_\xsa(U) \longrightarrow
  \mathcal{C}^{\infty,t}_{\xsa}(U)
  \overset{\bar\partial}{\longrightarrow}
  \mathcal{C}^{\infty,t}_{\xsa}(U)\longrightarrow 0 \ . \end{equation}

Now, we recall two results on the pull back of tempered holomorphic functions. We
refer to \cite{kashi} for the def\mbox{}inition of $\D_{X\to Y}$, for $f:X\to
Y$ a morphism of complex manifolds. For Lemma \ref{lem_pullback_ot},
see \cite[Lemma 7.4.7]{ast}.


\begin{lem}\label{lem_pullback_ot}
Let $f:X\to Y$ be a closed embedding of complex manifolds. There is a
natural isomorphism in $D^b(\rho_!\D_X)$
$$  \rho_!\D_{X\to Y}\underset{\rho_!f^{-1}\D_Y}{\overset{L}{\otimes}}f^{-1}\ot_Y\simeq\ot_X  \ . $$
\end{lem}

For Proposition \ref{prop_pullback_ot}, see \cite[Theorem 2.1]{ext1}.

\begin{prop}\label{prop_pullback_ot}
Let 
$f\in\mathcal O_\com(X)$, 
$U\in\opxsac$
such that $f|_{\overline U}$ is an injective map,
$h\in\mathcal O_X\bigr{f(U)}$. Then $h\in\ot_{\com_{sa}}\bigr{f(U)}$ if and
only if $h\circ f\in\otxsa(U)$. 
\end{prop}

We conclude this subsection by recalling the def\mbox{}inition of the sheaf
of holomorphic functions with moderate growth at the origin. We follow
\cite{malgr_birk}. Let $S^1$ be the unit circle,
$S^1\times\rea_{\geq0}$ the real blow-up at the origin of $\com^\times$.

For $\tau\in\rea$, $r\in\reagz$, $0<\ep<\pi$, the set
 $$ \csect\tau\ep r:=\big\{\rho e^{i\theta}\in\com^\times;\ \rho\in]0,r[,\ \theta\in]\tau-\ep,\tau+\ep[\big\} $$  
is called an {\em open sector centered at $\tau$ of amplitude $2\ep$ and radius
 $r$} or simply {\em an open sector}. Identifying $S^1$ with
$[0,2\pi[\subset\rea$, we will consider sectors centered at $\tau\in
S^1$. Further, with an abuse of language, we will say that an open sector $\csect\tau\ep r$ contains
$\theta\in\rea$ or $e^{i\theta}\in S^1$ if $\theta\in]\tau-\ep,\tau+\ep[\ (\mathrm{mod}\ 2\pi)$.

The sheaf on $S^1\times\rea_{\geq0}$ of holomorphic functions with
moderate growth at the origin, denoted $\mathcal A^{\leq0}$, is def\mbox{}ined as
follows. For $U$ an open set of $S^1\times\rea_{\geq0}$, set 
\begin{equation}\label{eq_moderate_growth}
\mathcal A^{\leq0}(U)    = \Big\{f
\parbox[t]{98mm}{$\in\mathcal O_\com(U\setminus(S^1\times\{0\}))$
  satisfying the following condition: for any $(e^{i\theta},0)\in U$
  there exist $C,M\in\reagz$ and an open sector $S\subset U$ containing
  $e^{i\theta}$ such that $\big|f(z)\big|< C|z|^{-M}$ for any $z\in
  S\Big\}$ \ .}
\end{equation}

Clearly $\mathcal A^{\leq0}$ is a sheaf on $S^1\times\rea_{\geq0}$.


In \cite{luca_micro}, the author def\mbox{}ines the functor $\nu_0^{sa}$ of
specialization at $0$ for the sheaves on the subanalytic site. One
has that, $\rho^{-1}\nu_0^{sa}(\ot_{\com_{sa}})\simeq\mathcal A^{\leq0}$.

\subsection{$\D$-modules on complex curves and $\Omega$-f\mbox{}iltered local systems}\label{dmod} 

In this subsection we recall some classical results on germs of $\D_X$-modules on
a complex analytic curve $X$. For a detailed and comprehensive
presentation we refer to \cite{malgr_birk}, \cite{kashi} and \cite{bv}.

Given a complex analytic curve $X$ and $x_0\in X$, we denote by $\mathcal O_X(*x_0)$
(resp. $\D_X(*x_0)$)
the sheaf on $X$ of holomorphic functions on $X\setminus\{x_0\}$,
meromorphic at $x_0$ (resp. the sheaf of rings of
dif\mbox{}ferential operators of f\mbox{}inite order with coef\mbox{}f\mbox{}icients in
$\mathcal O_X(*x_0)$). We set for short $\mathcal O_{*x_0}$ (resp. $\D_{*x_0}$) for
the stalk of $\mathcal O_X(*x_0)$ (resp. $\D_X(*x_0)$) at $x_0$. Further, we
denote by $\widehat{\mathcal O_{*x_0}}$ (resp. $\widehat{\D_{*x_0}}$) the
f\mbox{}ield of formal Laurent power series (resp. the ring of dif\mbox{}ferential
operators with coef\mbox{}f\mbox{}icients in $\widehat{\mathcal O_{*x_0}}$). The ring
$\mathcal O_{*x_0}$ comes equipped with a natural valuation
$v:\mathcal O_{*x_0}\to\integer\cup\{+\infty\}$.

By the choice of a local coordinate $z$ near $x_0$, we can suppose
that $X\subset\com$ is an open neighborhood of $x_0=0\in\com$. 

The category of holonomic $\D_{*0}$-modules, denoted
$\mod_h(\D_{*0})$, is equivalent to the category of local meromorphic
connections. 

For $\phi\in\zcz$, set $ \mc L^\phi:=\D_{*0}\exp(\phi)$.

For $l\in\integergz$, let $\mu_l:\com\to\com$, $z\mapsto z^l$. We
denote by $\mu_l^*$ the inverse image functor for
$\D_{*0}$-modules. 

Theorems \ref{thm_formal_LT} and \ref{thm_asymptotic_HT} below are
cornerstones in the theory of ordinary dif\-fe\-ren\-tial
equations. We refer to \cite{malgr_birk}, \cite{sabbah_iso} and
\cite{wasow}. 

\begin{thm}[Levelt-Turrittin's Formal Theorem]\label{thm_formal_LT}
Let $\mc M\in\mod_h(\D_{*0})$. There exist
$l\in\integergz$, a f\mbox{}inite set $\Sigma\subset\zcz$, a family of regular holonomic
  $\D_{*0}$-modules indexed by $\Sigma$, $\{\mc
  R_\phi\}_{\phi\in\Sigma}$, and an isomorphism in $\mod(\widehat{\D_{*0}})$
\begin{equation} \label{eq_formal_LT}
\mu_l^*\mc M\otimes\widehat{\mathcal O_{*0}}
\simeq
\ou{\phi\in\Sigma}{\oplus}\mc L^{\phi}\otimes\mc R_\phi\otimes\widehat{\mathcal O_{*0}}
\ .
\end{equation}
\end{thm}


In the litterature (for example \cite{malgr_birk}) the def\mbox{}inition of
the Katz invariant of $\mc M\in\mod_h(\D_{*0}) $ is given from the
Newton polygon of $\mc M$. For sake of simplicity we give an
equivalent def\mbox{}inition based on the isomorphism
\eqref{eq_formal_LT}. Clearly, the valuation $v$ induces a map, still
denoted $v$, $\zcz\to\integer\cup\{+\infty\}$.

\begin{df}
\begin{enumerate}
\item Let $\M\in\mod_h(\D_{*0})$. Suppose that \eqref{eq_formal_LT}
is satisf\mbox{}ied, $l$ is minimal and $\Sigma\neq\{0\}$. The
\emph{Katz invariant} of $\mc M$ is
$\underset{\phi\in\Sigma}{\max}\{\frac{-v(\phi)}l\}$.  
 If $\Sigma=\{0\}$ then the Katz invariant of $\mc
M$ is $0$. 

For $k\in\integergz$, we denote by $\mod_h(\D_{*0})_k$ the full abelian subcategory
of $\mod_h(\D_{*0})$ whose objects have Katz invariant strictly smaller than $k$. 
\item Let $\Sigma\subset\zcz$ be a f\mbox{}inite set and $\{\mc
  R_\phi\}_{\phi\in\Sigma}$ a family of regular holonomic
  $\D_{*0}$-modules indexed by $\Sigma$. A $\D_{*0}$-module isomorphic to
  $\ou{\phi\in\Sigma}{\oplus}\mc L^{\phi}\otimes\mc R_\phi$,  is said
  a \emph{good model}. 

  We denote by $\msf{GM}_k$ the full subcategory
  of $\mod_h(\D_{*0})_k$ whose objects are good models.  
\end{enumerate}
\end{df}


Roughly speaking, Theorem \ref{thm_asymptotic_HT} below says that the
formal isomorphism \eqref{eq_formal_LT} is analytic on
suf\mbox{}f\mbox{}iciently small open sectors.

Let
\begin{equation*}
 P:=\ou[m]{j=0}\sum a_j(z)\biggr{\frac d{dz}}^j
\  ,\end{equation*}
where $m\in\integer_{>0}$, $a_j\in\mathcal O_\com(X)$ and $a_m\not=0$.
 Denote by $\mc C^0_\com$ the sheaf of continuous functions
on $\com$. For $l\in\integergz$, $S$ an open sector,
$h\in\{1,\ldots,l\}$, let $\zeta_h:S\to\com$ be the $l$
dif\mbox{}ferent inverse functions to $z\mapsto z^l$ def\mbox{}ined on
$S$. 

\begin{thm}[Hukuhara-Turrittin's Asymptotic Theorem]\label{thm_asymptotic_HT}
There exist a finite set $\Sigma\subset\zcz$, $l,r_\phi\in\integergz$
($\phi\in\Sigma$) such that for any $\tau\in\rea$, there exist an open
sector $S$ containing $\tau$, $f_{\phi
  kh}\in\mathcal O_\com(S)\cap\mc{C}^0_\com\bigr{\overline{S}\setminus\{0\}}$
($\phi\in\Sigma,\,k=1,\ldots,r_\phi,\,h=1,\ldots,l$), satisfying 
\begin{enumerate}
\item $\bigc{f_{\phi kh}(z)\exp(\phi\circ\zeta_h(z));\ \phi\in\Sigma,\
    k=1,\ldots,r_\phi,\ h=1,\ldots,l}$ is a basis of the $\com$-vector 
  space of germs of holomorphic solutions of $Pu=0$ on $S$,
\item there exist $C,M\in\reagz$ such that, for any
  $z\in S$,
$$ C|z|^M\leq|f_{\phi kh}(z)|\leq\bigr{C|z|^M}^{-1} \  $$
\end{enumerate}
($\phi\in\Sigma,\ k=1,\ldots,r_\phi,\ h=1,\ldots,l$).
\end{thm}

As said above, the category $\mod_h(\D_{*0})$ is equivalent to the
category of local meromorphic connections. There already exists a
classif\mbox{}ication of local meromorphic connections by means of the
formal structure and the Stokes phenomena. We recall it in the language of
$\Omega$-graded and $\Omega$-f\mbox{}iltered local systems. We refer
to \cite{dmr}, \cite{malgr_class}, \cite{malgr_birk} and \cite{bv}.

Let $\Omega^1$ be the sheaf on $\com^\times$ of $1$-forms with
holomorphic coef\mbox{}f\mbox{}icients. Let $\Omega$ be the
local system on $S^1$ def\mbox{}ined as follows: for $\theta\in S^1$,
choose a determination of $z\mapsto z^{1/l}$ near $\theta$ and set 

\begin{multline*}
\Omega_\theta:=\Big\{\ou[n]{j=1}{\sum}a_jz^{-(\frac jl+1)}dz\in\Omega^1(S); \\
a_j\in\com,\ l\in\integergz, \ S \textrm{ an open sector containing }\theta \Big\}  \ .
\end{multline*}

The $\com$-vector space $\Omega_\theta$ can be endowed with a partial
order. For $\alpha,\beta\in \Omega_\theta$, we write
$\beta\prec_\theta\alpha$, if $\exp\bigr{\int\beta-\alpha}\in\mc
A^{\leq0}_{(\theta,0)}$. 

\begin{df}
An $\Omega$-\emph{graded local system} on an open subset $U$ of $S^1$ is the
data of a f\mbox{}inite rank local system $\mc V$ on $U$ such that for any
$\theta\in U$, $\mc V_\theta$ is an $\Omega_\theta$-graded
$\com$-vector space, $\mc
V_\theta\simeq\underset{\alpha\in\Omega_\theta}{\oplus}\mc
V_{\theta,\alpha}$, satisfying the following condition. If $\theta$and
$\theta'$ are in the same connected component of $U$, then the grading of
$\mc V_{\theta'}$ is induced by the grading of $\mc V_{\theta}$ by
analytic continuation. 

A morphism between two $\Omega$-graded local systems $\mc V_1,\mc V_2$
is a morphism of local systems $f:\mc V_1\to\mc V_2$ such that, for
any $\theta\in S^1$, $f_\theta:V_{1\theta}\to\mc V_{2\theta}$ is a
morphism of $\Omega_\theta$-graded $\com$-vector spaces. 
\end{df}

The analytic continuation induces a monodromy action of $\integer$ on
$\Omega_\theta$ that we denote by $\mu_\Omega$. Given an
$\Omega$-graded local system $\mc V$, the monodromy action of
$\integer$ on $\mc V_\theta$ induces an isomorphism $\mc
V_{\theta,\alpha}\simeq\mc V_{\theta,\mu_\Omega\,\alpha}$. Further,
for $\theta\in S^1$, the functor $\mc V\mapsto\mc V_\theta$ is an
equivalence between the category of $\Omega$-graded local systems and
the category of f\mbox{}inite dimensional $\Omega_\theta$-graded
$\com$-vector spaces equipped with a $\integer$-action compatible with
$\mu_\Omega$.

Let us now def\mbox{}ine $\Omega$-f\mbox{}iltered local systems. Roughly speaking,
an $\Omega$-f\mbox{}iltered local system is obtained by gluing
$\Omega$-graded local systems with maps which preserve the
$\Omega$-f\mbox{}iltration induced by the partial order on $\Omega_\theta$.

\begin{df}
An $\Omega$-\emph{f\mbox{}iltered local system} is a f\mbox{}inite rank local system $\mc
V$ obtained by the following gluing data. 
\begin{enumerate}
\item An open covering $\{U_j\}_{j\in J}$ of $S^1$ and $\Omega$-graded
  local systems $\mc V_j$ on $U_j$, $\mc
  V_{j\theta}\simeq\underset{\alpha\in\Omega_\theta}{\oplus}\mc
  V_{j\theta,\alpha}$.  
\item A family of isomorphisms of local systems $\phi_{jk}:\mc V_{k}|_{U_j\cap
    U_k}\overset\sim\lra\mc V_{j}|_{U_j\cap U_k}$ such that, for any
  $\theta\in U_j\cap U_k$,  
\begin{equation}\label{eq_gluing_graded}  \phi_{jk\theta}\big|_{\mc
  V_{k\theta,\alpha}}:\mc
  V_{k\theta,\alpha}\lra\mc
  V_{j\theta,\alpha}\oplus\Bigr{\underset{\beta\prec_\theta\alpha}{\oplus}\mc V_{j\theta,\beta}}   \ . \end{equation} 
\end{enumerate} 
\end{df}

Remark that composing the morphism \eqref{eq_gluing_graded} with the
projection on $\mc V_{j\theta,\alpha}$, one obtains an invertible morphism $\mc V_{k\theta,\alpha}\to\mc V_{j\theta,\alpha}$. Thus one def\mbox{}ines invertible morphisms of $\Omega$-graded $\com$-vector spaces $\mc V_{k\theta}\to\mc V_{j\theta}$.  In this way one def\mbox{}ines the functor
$gr$ from the category of $\Omega$-f\mbox{}iltered local systems to the
category of $\Omega$-graded local systems. 


Thanks to Theorem \ref{thm_asymptotic_HT}, for $\mc
M\in\mod_h(\D_{*0})$, the solutions of $\mc M$ can be naturally
endowed with a structure of $\Omega$-f\mbox{}iltered local system. Thus we
have a functor, $\ms S^\Omega$, from the category $\mod_h(\D_{*0})$ to
the category of $\Omega$-f\mbox{}iltered local systems.  

\begin{df}\label{df_formal_invariants}
Let $\mc M\in\mod_h(\D_{*0})$, the elements of the set
$$ \Omega(\mc M):=\Bigc{\phi\in\zlczl\textrm{ such that } gr_{\frac{d\phi}{dz} dz}(\ms S^\Omega(\mc M)_\theta)\neq 0} $$
are called the \emph{determinant polynomials} of $\mc M$. Further we set
$$ r_{\phi,\mc M}:=\dim\ gr_{\frac{d\phi}{dz} dz}\bigr{\ms S^\Omega(\mc M)_\theta} \ .$$
\end{df}

For the proof of Theorem \ref{thm_I-filt} below we refer to \cite{malgr_birk} or \cite{bv}.

\begin{thm}\label{thm_I-filt}
\begin{enumerate}
\item The functor $\ms S^\Omega$ from the category $\mod_h(\D_{*0})$
  to the category of $\Omega$-f\mbox{}iltered local systems is an equivalence
  of categories. 
\item The functor $\ms S^\Omega$ restricted to the category
  $\mod_h(\widehat{\D_{*0}})$ is an equivalence of categories onto the
  category of $\Omega$-graded local systems. 
\end{enumerate}
\end{thm}


\section{Tempered growth of exponential functions}\label{SEC_EXP}

This section is subdivided as follows. In the f\mbox{}irst part we study the
family of sets where a function of the form $\exp(\phi)$,
$\phi\in\zcz$, is tempered. In the second part we use the results of
the f\mbox{}irst part in order to prove that such family determines $\phi$ up
to a multiplicative positive constant. Throughout this section
$X=\com$.

\subsection{Sets where exponential functions have
  tempered growth}

For $\phi\in\zcz\setminus\{0\}$ and $A\in\reagz$, set 
\begin{equation}\label{eq_UphiA}  
U_{\phi,A}:=\big\{z\in\com^\times;\
\re\big(\phi(z)\big)<A\big\} \ ,  
\end{equation}
further set $U_{0,A}:=\com$.

First we state and prove the analogue of a result of \cite{arx}.

\begin{prop}\label{prop_generalization_arx}
  Let $\phi\in\zcz$ and $U\in\opxsac$ with $U\neq\varnothing$.
  The conditions below are equivalent.
  \begin{enumerate}
  \item $\exp\big(\phi\big)\in\otxsa(U)$.
  \item There exists $A\in\reagz$ such that $U\subset U_{\phi,A}$.
  \end{enumerate}
\end{prop}

Before proving Proposition \ref{prop_generalization_arx}, we need the
following

\begin{lem}[\cite{arx}]\label{lem_arx}
  Let $W\neq\varnothing$ be an open subanalytic subset of $\mathbb
  P^1(\com)$, $\infty\notin W$. The following conditions are
  equivalent.
  \begin{enumerate}
  \item There exists $A\in\reagz$ such that $\re z<A$, for any $z\in
    W$.
  \item The function $\exp(z)$ has polynomial growth on any
    semi-analytic arc $\Gamma\subset W$ with an endpoint at $\infty$.
    That is, for any semi-analytic arc $\Gamma\subset W$ with an
    endpoint at $\infty$, there exist $M,C\in\reagz$ such that, for
    any $z\in\Gamma$,
    \begin{equation}\label{eq_lem_arx}
      \big|\exp(z)\big|\leq C\bigr{1+|z|^2}^M \ .
    \end{equation}
  \end{enumerate}
\end{lem}

\emph{Proof.}  Clearly, \emph{(i)}$\Rightarrow$\emph{(ii)}.

Let us prove \emph{(ii)}$\Rightarrow$\emph{(i)}. Set $z:=x+iy$ and
suppose that $x$ is not bounded on $W$. 
There exist $\ep,L\in\reagz$ and
a real analytic map 
\begin{eqnarray*}
  \gamma:[0,\ep[    &  \longrightarrow  &  \mathbb P^1(\com)  \\
  t                 &  \longmapsto      &  (x(t),y(t)) \ ,
\end{eqnarray*} 
such that
$\gamma(0)=\infty$, $\gamma\big(]0,\ep[\big)\subset W$ and
$x\bigr{]0,\ep[}=]L,+\infty[$. Since $\gamma$ is analytic, there exist $q\in\rat$,
$c\in\rea$ and $\mu\in\reagz$ such that, for any $t\in]0,\ep[$,
$$\gamma(t)=\Bigr{x(t)\ ,\ c\,x(t)^q+O\bigr{x(t)^{q-\mu}}} \ .$$
Now, if \eqref{eq_lem_arx} is satisf\mbox{}ied, then $\exp(x)$ has polynomial
growth in a neighborhood of $+\infty$, which gives a contradiction.
\proofend

\emph{Proof of Proposition \ref{prop_generalization_arx}.}

Clearly, $(ii)\Rightarrow(i)$.

$(i)\Rightarrow(ii)$. 
The result is obvious if $\phi=0$. Otherwise, we distinguish two cases.

\begin{underline}{Case 1}\end{underline}:
suppose $\phi(z)=\frac1z$.

Suppose that for any $A\in\reagz$, there exists $z_A\in U$ such
that $\re\big(\frac1{z_A}\big)>A$. Then, by Lemma \ref{lem_arx}, there exists a
semi-analytic arc with an endpoint at $0$, $\Gamma\subset U$, such
that $\exp\big(\re \frac1{z}\big)$ has not polynomial growth on
$\Gamma$. That is, for any
$M,C\in\reagz$, there exists $z_{M,C}\in\Gamma$ satisfying
$$ \exp\big(\re\frac1{z_{M,C}}\big)\geq\frac C{|z_{M,C}|^M} \ .$$


Apply Proposition \ref{prop_reg_dist} with $X=\overline{\Gamma}$
and $Y=\partial U$. There exist an open neighborhood $V$
of $0$, $c,r\in\reagz$ such that, for any $z\in\Gamma\cap V$,
$$  |z|\leq c\,\dist z{\partial U}^r  \ .$$ 

Hence,
$$ \exp\big(\re\frac1{z_{M,C}}\big)\geq\frac {c^{-1}C}{\dist
  {z_{M,C}}{\partial U}^{rM}}  \ .$$

It follows that $\exp\big(\frac1z\big)$ is not tempered on $U$.



\begin{underline}{Case 2}\end{underline}:
suppose $ \phi(z)=\ou[n]{j=1}\sum\frac{a_j}{z^j} $, with
$n\in\integergz$ and $a_n\neq0$. 







Let
\begin{eqnarray*}
\eta(z)  & :=  &  \biggr{\ou[n]{j=1}\sum\frac{a_j}{z^j}}^{-1} \\
  &  =  &  \frac{z^n}{\ou[n]{j=1}\sum a_j z^{n-j}}  \ .
\end{eqnarray*}

There exists a neighborhood $W\subset\com$ of $0$ such that
$\eta\in\mathcal O_\com(W)$. It is well known that a non-constant holomorphic function
is locally the composition of a holomorphic isomorphism and a positive
integer power of $z$. Since it is suf\mbox{}f\mbox{}icient to prove the result in a
neighborhood of $0$ and up to f\mbox{}inite coverings, we can suppose that
$U\subset W$ and $\eta|_{\overline U}$ is injective.

Consider the following commutative triangle,
$$ 
\xymatrix{
U\ar[r]^\eta\ar[rd]_{\exp(\phi)} & \eta(U)\ar[d]^{\exp(\frac1\zeta)}  \\  &  \com \ .
} 
$$

Using Proposition \ref{prop_pullback_ot} and Case 1, we have that

\begin{eqnarray*}
\exp(\phi)\in\ot(U)  & \Leftrightarrow  &  \exp\biggr{\frac1\zeta}\in\ot(\eta(U))  \\
  &  \Leftrightarrow  &  \eta(U)\subset U_{\frac1\zeta,A} \quad
  \textrm{for some } A\in\reagz\\
  &  \Leftrightarrow  &  U\subset U_{\phi,A}  \ \quad\quad
  \textrm{for some } A\in\reagz \ .
\end{eqnarray*}

\proofend

\begin{cor}\label{cor_generalization_arx_ram}
  Let $\phi\in\zcz$. Let $S$ be an open sector of amplitude $2\pi$,
  $U\in\opxsac$, $\varnothing\neq U\subset S$, $\zeta:S\to\com$, an inverse of
  $z\mapsto z^{l}$. Then $\exp\big(\phi\circ\zeta\big)\in\otxsa(U)$ if
  and only if there exists $A\in\reagz$ such that $\zeta(U)\subset
  U_{\phi,A}$.

  In particular, setting
$$  U_{\phi\circ\zeta,A}:=\Bigc{z\in S;\ \re\bigr{\phi\circ\zeta(z)}<A}   \ , $$ 
one has that $\exp(\phi\circ\zeta)\in\otxsa(U)$ if and only if there
exists $A\in\reagz$ such that $U\subset U_{\phi\circ\zeta,A}$.
\end{cor}

\emph{Proof.} 
Let $\mu_l(z):=z^l$. Consider the following commutative
diagram,
$$ 
\xymatrix{
\zeta(U)\ar[r]^{\mu_l}\ar[rd]_{\exp(\phi)} &
U\ar[d]^{\exp(\phi\circ\zeta)}\\
  & \com \ .
}
 $$

By Proposition \ref{prop_pullback_ot}, we have that
$\exp(\phi\circ\zeta)\in\ot_\xsa(U)$ if and only if
$\exp(\phi)\in\ot_\xsa(\zeta(U))$. Then the conclusion follows by
Proposition \ref{prop_generalization_arx}. 

\proofend

Now we are going to introduce a class of subanalytic sets which plays
an important role in what follows. 

\begin{df} \label{df_tau_concentrated}
For $\tau\in\rea$ we say that $U\in\opxsac$ is \emph{concentrated
  along $\tau$} if $U\neq\varnothing$ is connected, $0\in\partial U$ and, for any
 open sector $S$ containing $\tau$, there exists an open neighborhood $W\subset\com$ of
  $0$ such that $ U\cap W\subset S$.
\end{df}

Lemma \ref{lem_tau_concentrated_thru_holomorphic} below follows easily
from the well known fact that a non-constant holomorphic function
is locally the composition of a holomorphic isomorphism and a positive
integer power of $z$.

\begin{lem}\label{lem_tau_concentrated_thru_holomorphic}
Let $W\subset\com$ be an open neighborhood of $0$, $f\in\mathcal
O(W)$. Suppose that $f$ has a zero of order
  $l\in\integergz$ at $0$. There exists $\tau_f\in\rea$, depending
  only on the argument of $f^{(l)}(0)$, satisfying the
  following conditions.
\begin{enumerate}
\item \sloppy For any $\tau\in\rea$, $U\in\opxsac$ concentrated along $\tau$,
  there exists an open neighborhood $W'$ of $0$, $\overline{W'}\subset W$, such that $f|_{\overline {U\cap W'}}$ is injective and
  $f(U\cap W')$ is concentrated along $l(\tau+\tau_f)$.
\item For any $\tau\in\rea$, $V\in\opxsac$ concentrated along $\tau$,
  there exist an open neighborhood $W'$ of $0$, $V_0,\ldots,V_{l-1}\in\opxsac$, such that 
 such that $V_j$ is concentrated along $\frac\tau l-\tau_f+j{\frac{2\pi}l}$, 
 and $f(V_j)=V\cap W'$ 
($j=0,\ldots,l-1$).
\end{enumerate}
\end{lem}


Proposition \ref{prop_+-temp} below will play a fundamental role in the next
subsection. 

\begin{prop}\label{prop_+-temp}
  Let $n\in\integer_{\geq0}$, $\tau_0\in\rea$. There exists
  $\tau\in\rea$ such that, for any $\phi=\frac{\rho
    e^{i\tau_0}}{z^n}+\widetilde\phi\in\zcz$ ($\rho\in\reagz,\
  -v(\widetilde\phi)<n$), there exist $U_0,\ldots,U_{2n-1}\in\opxsac$
  satisfying
  \begin{enumerate}
  \item $U_j$ is concentrated along $\tau+
    j\frac\pi n$,
  \item $\exp\bigr{\phi},\exp\bigr{-\phi}\in\otxsa(U_j)$,
  \end{enumerate}
  ($j=0,\ldots,2n-1$).
\end{prop}

\emph{Proof.} The result is obvious if $\phi=0$. Otherwise we distinguish three cases.
\begin{underline}{Case 1}\end{underline}: suppose $\phi(z)=\frac1z$.

Recall \eqref{eq_UphiA}. For $A\in\reagz$, one checks easily that the
set $U_{\frac1z,A}$ (resp. $U_{-\frac1z,A}$) is the complementary of
the closed ball of center $\big({\frac1{2A}},0\big)$ (resp.
$\big(-{\frac1{2A}},0\big)$) and radius ${\frac1{2A}}$.

Set
\begin{eqnarray*}
U_1  &  :=  &  \big\{(x,y)\in\rea^2; \ |x|<1, \ \sqrt{|x|-x^2}<y<1 \big\} \ ,\\
U_2  &  :=  &  \big\{(x,y)\in\rea^2; \ |x|<1, \ -1<y<-\sqrt{|x|-x^2} \big\} \
.
\end{eqnarray*}

It is easy to see that $U_1$ (resp. $U_2$) is concentrated
along ${\frac\pi2}$ (resp. ${\frac{3\pi}2}$) and $U_1\cup U_2\subset
U_{\frac1z,1}\cap U_{-\frac1z,1}$. Hence, by Proposition
\ref{prop_generalization_arx}, 
\begin{equation}\label{eq_ot_u1u2} 
\exp(1/z),\ \exp(-1/z)\in\otxsa(U_j)\quad\quad(j=1,2) \ .
\end{equation}

\begin{underline}{Case 2}\end{underline}: suppose that
$\phi(z)=\frac1{z^m}$, for $m\in\integergz$. Let $\mu_m:\com\to\com$,
$\mu_m(z)=z^m$. Consider the commutative triangle
$$ \xymatrix{
\com^\times\ar[rd]_{\exp{(1/z^m)}}\ar[r]^{\mu_m}  &  \com^\times\ar[d]^{\exp(1/z)}\\ & \com \ .
} $$

Consider $U_1,U_2$ as in Case 1. Applying Lemma
\ref{lem_tau_concentrated_thru_holomorphic} \emph{(ii)} with
$f=\mu_m$, we obtain that there exist
$V_{1,0},\ldots,V_{1,m-1}\in\opxsac$ (resp.
$V_{2,0},\ldots,V_{2,m-1}\in\opxsac$) such that
\begin{enumerate} 
\item $V_{1,j}$ (resp. $V_{2,j}$) is concentrated along
  $\frac\pi{2m}+j\frac{2\pi}m$ (resp. $\frac{3\pi}{2m}+j\frac{2\pi}m$), 
\item $\mu_m(V_{k,j})=U_k$\hspace{0.2cm}, 
\end{enumerate}
 ($j=0,\ldots,m-1,\ k=1,2$).

Clearly, $\mu_m|_{\overline{V_{k,j}}}$ is injective. 
By Proposition \ref{prop_pullback_ot}, we have that 
$$ \exp\bigr{1/z^m}\in\otxsa(V_{k,j}) $$
$$ \Bigr{\text{resp. } \exp\bigr{-1/z^m}\in\otxsa(V_{k,j})} $$
if and only if
$$ \exp\bigr{1/z}\in\otxsa\bigr{\mu_m(V_{k,j})}=\otxsa(U_k) $$
$$ \Bigr{\text{resp. } \exp\bigr{-1/z}\in\otxsa\bigr{\mu_m(V_{k,j})}=\otxsa(U_k)} $$
$(j=0,\ldots,m-1,k=1,2)$. The conclusion follows.

\begin{underline}{Case 3}\end{underline}: suppose that  
$$ \phi(z)=\sum_{j=1}^n \frac{a_j}{z^j}\in\zcz \ ,$$
for $n\in\integergz$, $a_j\in\com$ ($j=1,\ldots,n$) and $a_n\neq0$.

First, we recall the implicit function theorem for convergent power
series. We denote by $\com\{x\}$ (resp. $\com\{x,y\}$) the ring of
convergent power series in $x$ (resp. $x,y$). We refer to \cite[Theorem 8.6.1, p.166]{che} for the proof. 

\begin{thm}\label{thm_puiseux_convergent_series}
Let $F\in\com\{x,y\}$ be such that $F(0,0)=0$. There exists
$\eta(x)\in\cup_{l\in\integergz}x^{1/l}\com\{x^{1/l}\}$ such that $F(x,\eta(x))=0$.
\end{thm}

Consider
\begin{eqnarray}\label{eq_poly_phi_zeta}
F(z,\eta)  &  :=  &  -\eta^{n}+z^n\sum_{j=1}^{n}a_j\eta^{n-j}\
\notag \\
  &  =  &
  -\eta^n+z^na_1\eta^{n-1}+\ldots+z^na_{n-1}\eta+z^na_n \ .
\end{eqnarray}
By Theorem \ref{thm_puiseux_convergent_series}, there exist $l\in\integergz$, 
$\eta(z)\in z^{1/l}\com\{z^{1/l}\}$ such that
$F\big(z,\eta(z)\big)=0$. Since $a_{n}\neq0$, we have that
$\eta(z)\neq0$, for $z\neq0$. It follows that $\eta(z)\in z^{1/l}\com\{z^{1/l}\}$ satisf\mbox{}ies 
\begin{equation}\label{eq_p_circ_phi} 
\phi\big(\eta(z)\big)=\sum_{j=1}^{n}\frac{a_j}
  {\eta(z)^{j}}={\frac1{z^n}} \ . 
\end{equation}

Further, substituting $\eta(z)$ in \eqref{eq_poly_phi_zeta}, one checks
that $l=1$ and 
 $\eta(z)=z\sigma(z) $, 
for $\sigma$ an invertible element of $\com\{z\}$ such that 
$\arg\bigr{\sigma(0)}=\frac{\arg(a_n)}{n}$. In particular, 
 there exists an open neighborhood $W\subset\com$ of the 
origin such that $\eta\in\mathcal O_\com(W)$. 

Now, by Case 2, there exist $V_{k,j}\subset W$ $(j=0,\ldots,n-1,\ k=1,2)$
such that $V_{1,j}$ (resp. $V_{2,j}$) is concentrated along
  $\frac\pi{2n}+j\frac{2\pi}{n}$
  (resp. $\frac{3\pi}{2n}+j\frac{2\pi}{n}$) and

\begin{equation}\label{eq_exp1/znl}
  \exp\Bigr{\frac1{z^{n}}},\exp\Bigr{-\frac1{z^{n}}}\in\otxsa(V_{k,j}) 
\end{equation}
$(j=0,\ldots,n-1,\ k=1,2)$.


As $\eta$ has a zero of order $1$ at $0$, by Lemma \ref{lem_tau_concentrated_thru_holomorphic} \emph{(i)}, there
exists $\tau_\eta\in\rea$, depending only on
$\arg\bigr{\eta(0)}=\frac {\arg(a_n)}n$, such that, up to shrinking
$W$,
\begin{enumerate}
\item $\eta|_{\overline{V_{k,j}}}$ is injective and
\item $\eta\bigr{V_{1,j}}$ (resp. $\eta\bigr{V_{2,j}}$) is
  concentrated along $\tau_\eta+\frac\pi{2n}+j\frac{2\pi}{n}$
  (resp. $\tau_\eta+\frac{3\pi}{2n}+j\frac{2\pi}{n}$),
\end{enumerate} 
$(j=0,\ldots,n-1)$.

Consider the commutative triangle
$$ 
\xymatrix{
  W\setminus\{0\}\ar[r]^\eta\ar[rd]_{\exp\bigr{\frac1{z^{n}}}} &
  \com^\times\ar[d]^{\exp(\phi(z))}
  \\
  & \com \ .  }
$$

By Proposition \ref{prop_pullback_ot}, we have that 
$$ \exp\bigr{\phi(z)}\in\otxsa\bigr{\eta(V_{k,j})} $$
$$ \Bigr{\text{resp. } \exp\bigr{-\phi(z)}\in\otxsa\bigr{\eta(V_{k,j})}} $$
if and only if
$$ \exp\bigr{\phi\circ\eta(z)}=\exp\bigr{1/z^{n}}\in\otxsa\bigr{V_{k,j}} $$
$$ \Bigr{\text{resp. } \exp\bigr{-\phi\circ\eta(z)}=\exp\bigr{-1/z^{n}}\in\otxsa\bigr{V_{k,j}}}$$
$(j=0,\ldots,n-1,k=1,2)$. The conclusion follows from \eqref{eq_exp1/znl}.

\nopagebreak
\proofend

\begin{rem}
{\rm Recall the def\mbox{}inition given in the end of Subsection
  \ref{SUBSEC_OT} of the sheaf $\mathcal A^{\leq0}$ def\mbox{}ined on
  $S^1\times\rea_{\geq0}$, considered as the real blow-up at $0$ of
  $\com^\times$. Let $\tau\in\rea$, $U\in\opxsac$ concentrated along
  $\tau$, the set $(\tau,0)\cup U\subset S^1\times\rea_{\geq0}$ is not
  open.  Further if $\exp(\phi)\in\mathcal A^{\leq0}_{(\tau,0)}$ then
  $\exp(-\phi)\notin\mathcal A^{\leq0}_{(\tau,0)}$.}
\end{rem}

We conclude this subsection with an easy lemma which will be useful in
the next subsection. First, let us introduce some notation.

Given $\phi\in\zcz$, $\phi={\eta e^{i\tau}\over z^n}+\tilde\phi$, for
$\eta\in\reagz,n\in\integergz,\tau\in\rea$ and $\tilde\phi\in\zcz$, $-v(\tilde\phi)<n$, set
$$ I_\phi:=\Big\{\theta\in[0,2\pi];\ \cos(\tau-n\theta)<0 \Big\} \
. $$ 
In other words, $I_\phi$ is the support of $\exp(\phi)$ as a
section of $\mc A^{\leq0}\big|_{S^1\times\{0\}}$.

Recall the def\mbox{}inition of the sets $U_{\phi,A}$ given in
\eqref{eq_UphiA}. 

\begin{lem}\label{lem_composed}
\begin{enumerate}
\item Let $\phi_1,\phi_2\in\zcz$, $\phi_j:={\eta_j e^{i\tau_j}\over
    z^{n_j}}+\tilde\phi_j$, for $\eta_j\in\reagz, n_j\in\integergz,\tau_j\in\rea$
  and $\tilde\phi_j\in\zcz$, $-v(\tilde\phi_j)<n_j$.

  If $n_1\neq n_2$ or $\tau_1\neq\tau_2$, then
  $I_{\phi_1}\setminus\overline I_{\phi_2}\neq\varnothing $ and
  $I_{\phi_2}\setminus\overline I_{\phi_1}\neq\varnothing $.
\item Let $\theta_0\in[0,2\pi[$ and $\phi\in\zcz\setminus\{0\}$.  If
  $\theta_0\in I_{\phi}$, then there exists an open sector $S$
  containing $\theta$ such that, for any $A\in\reagz$, $ S\subset U_{\phi,A} $.
In particular, for any $U\in\opxsac$ concentrated along $\theta_0$,
$\exp(\phi)\in\otxsa(U)$. 
\item Let $\theta_0\in[0,2\pi[$ and $\phi\in\zcz\setminus\{0\}$. 
If $\theta_0\notin\overline I_{\phi}$, then there exists an open sector $S$ containing $\theta$
such that, for any $A\in\reagz$, $ S\subset X\setminus U_{\phi,A} $. 
In particular, for any $U\in\opxsac$ concentrated along $\theta_0$,
$\exp(\phi)\notin\otxsa(U)$.
\end{enumerate}
\end{lem}

\emph{Proof.}

The result follows from some easy computations.

\proofend

\subsection{Comparison between growth of exponential functions}

In this subsection we are going to use the results of the previous
subsection in order to prove that if $\phi_1,\phi_2\in\zcz$ and, for any
$\lambda\in\reagz$, $\phi_1\neq\lambda\phi_2$, then the families
$\{U_{\phi_1,A}\}_{A\in\reagz}$ and $\{U_{\phi_2,A}\}_{A\in\reagz}$
are not cof\mbox{}inal.

The main result of this subsection is Proposition \ref{prop_exp_case}
below.

\begin{prop}\label{prop_exp_case}
Let $\phi_1,\phi_2\in\zcz\setminus\{0\}$. 
\begin{enumerate}
\item Suppose that there exists $\lambda\in\reagz$ such that
  $\phi_1=\lambda\phi_2$. Then, for any $A\in\reagz$,
  $U_{\phi_1,A}=U_{\phi_2,\frac A\lambda}$. In particular, for any
  $U\in\opxsac$, $\exp(\phi_1)\in\ot_\xsa(U)$ if and only if
  $\exp(\phi_2)\in\ot_\xsa(U)$. 
 \item Suppose that, for any $\lambda\in\reagz$,
$\phi_1\neq\lambda\phi_2$. Then for any open sector $S$ of amplitude
$>\frac{2\pi}{\max\{-v(\phi_1),-v(\phi_2),2\}}$ at least one of the two
conditions below is satisf\mbox{}ied (resp. for any open neighborhood $S$ of $0$, both conditions below are satisf\mbox{}ied).
\begin{enumerate}
\item There exists $U\in\opxsac$, $U\subset S$ such that
  $\exp(\phi_1)\in\ot_\xsa(U)$ and $\exp(\phi_2)\notin\ot_\xsa(U)$.
\item There exists $V\in\opxsac$, $V\subset S$ such that
  $\exp(\phi_1)\notin\ot_\xsa(V)$ and $\exp(\phi_2)\in\ot_\xsa(V)$.
\end{enumerate}

\end{enumerate}
\end{prop} 

\emph{Proof
.} 

\emph{(i)} Obvious.

\emph{(ii)} For $S$ an open sector, set $\tilde
S:=\bigc{\theta\in[0,2\pi[; \exists\ r>0\ re^{i\theta}\in S}$.

Let
$$ 
\phi_1(z)  :=  {\eta_1 e^{i\tau_1}\over z^{n_1}}+\tilde \phi_1(z)
\qquad\textrm{ and }\qquad  \phi_2(z)  :=  {\eta_2 e^{i\tau_2}\over z^{n_2}}+\tilde \phi_2(z)
\ ,
$$
for $\eta_j\in\reagz$, $\tau_j\in[0,2\pi[$ and $\tilde
\phi_j(z)\in\zcz$, $-v(\tilde \phi_j)<n_j$ ($j=1,2$).

Suppose that $n_1\neq n_2$ or $\tau_{1}\neq\tau_{2}$.

By Lemma \ref{lem_composed} \emph{(i)}, $I_{\phi_1}\setminus\overline
I_{\phi_2}\neq\varnothing$ and $I_{\phi_2}\setminus\overline I_{\phi_1}\neq\varnothing$.

Let $S$ be an open neighborhood of $0$, $\theta_1\in
I_{\phi_1}\setminus\overline I_{\phi_2}$ and $\theta_2\in
I_{\phi_2}\setminus\overline I_{\phi_1}$. There exist $U_k\in\opxsac$
concentrated along $\theta_k$ such that $U_k\subset S$ ($k=1,2$). The
result follows by Lemma \ref{lem_composed} \emph{(ii),(iii)}.

Suppose that $S$ is an open sector of amplitude
$>\frac{2\pi}{\max\{n_1,n_2,2\}}$. Then, there exists $\theta\in\tilde
S$ such that either $\theta\in I_{\phi_1}\setminus\overline
I_{\phi_2}$ or $\theta\in I_{\phi_2}\setminus\overline I_{\phi_1}$.
Since $\theta\in\tilde S$, there exists $U\in\opxsac$ concentrated
along $\theta$ such that $U\subset S$. The conclusion follows by Lemma
\ref{lem_composed} \emph{(ii),(iii)}.

 Now suppose that $n_1=n_2=n$ and $\tau_{1}=\tau_{2}$. That is,

$$ 
\phi_1(z)=\frac{\eta_1e^{i\tau_{1}}}{z^n}+\tilde \phi_1(z)
\qquad\textrm{ and }\qquad  \phi_2(z)  =\frac{\eta_2e^{i\tau_{1}}}{
  z^n}+\tilde \phi_2(z) 
\ .
$$

Since, for any $\lambda\in\reagz$, $\phi_1\neq\lambda \phi_2$, we have
that $n\geq2$. 

Set $\psi_{21}:=\phi_2-\frac{\eta_2}{\eta_1}\phi_1$ and
$\psi_{12}:=\phi_1-\frac{\eta_1}{\eta_2}\phi_2$. Since
$\psi_{21}\neq0$ and
$\psi_{21}=-\frac{\eta_2}{\eta_1}\bigr{\phi_1-\frac{\eta_1}{\eta_2}\phi_2}=-\frac{\eta_2}{\eta_1}\psi_{12}$,
then $I_{\psi_{21}}=I_{-\psi_{12}}$.

By Proposition \ref{prop_+-temp}, there exist
$\tau\in\rea$ and $U_0,\ldots,U_{2n-1},V_0,\ldots,V_{2n-1}\in\opxsac$ satisfying the conditions
\begin{enumerate}
\item $U_j, V_j$ are concentrated along $\tau+j\frac\pi n$,
\item $\exp(\phi_1),\exp(-\phi_1)\in\otxsa(U_j)$ and $\exp(\phi_2),\exp(-\phi_2)\in\otxsa(V_j)$,
\end{enumerate} 
($j=0,\ldots,2n-1$).

Since $-v(\psi_{12})=-v(\psi_{21})<n$, if $S$ is an open sector of
amplitude $>\frac{2\pi}n$, there exists $j'\in\{0,\ldots,2n-1\}$ such
that $\tau+j'\frac\pi n\in\tilde S$ and either $\tau+j'\frac\pi
n\notin\overline I_{\psi_{12}}$ or $\tau+j'\frac\pi n\notin\overline
I_{\psi_{21}}$.

More generically, $\{\tau+j\frac\pi n;\
j\in0,\ldots,2n-1\}\not\subset\overline I_{\psi_{12}}$ and
$\{\tau+j\frac\pi n;\ j\in0,\ldots,2n-1\}\not\subset\overline
I_{\psi_{21}}$.


Let us consider the case $\tau+j'\frac\pi n\notin\overline
I_{\psi_{12}}$. Since $V_{j'}$ is concentrated along $\tau+j'\frac\pi
n$, Lemma \ref{lem_composed} \emph{(iii)} implies
\begin{equation}\label{eq_psi12}
\exp(\psi_{12})\notin\ot_\xsa(V_{j'}) \ .
\end{equation}
Suppose now that $\exp(\phi_1)\in\ot_\xsa(V_{j'})$. Since
$\exp(-\phi_2)$, $\exp(\phi_2)\in\ot_\xsa(V_{j'})$ and the product of
tempered functions is tempered, we have that
$\exp(\phi_1-\frac{\eta_1}{\eta_2}\phi_2)\in\ot_\xsa(V_{j'})$, which
contradicts \eqref{eq_psi12}. Hence
$\exp(\phi_1)\notin\ot_\xsa(V_{j'})$.

Let us consider the case $\tau+j'\frac\pi n\notin\overline
I_{\psi_{21}}$. Since $U_{j'}$ is concentrated along $\tau+j'\frac\pi
n$, Lemma \ref{lem_composed} \emph{(iii)} implies
\begin{equation}\label{eq_psi21}
\exp(\psi_{21})\notin\ot_\xsa(U_{j'}) \ .
\end{equation}
Suppose now that $\exp(\phi_2)\in\ot_\xsa(U_{j'})$. Since
$\exp(\phi_1)$, $\exp(-\phi_1)\in\ot_\xsa(U_{j'})$ and the product of
tempered functions is tempered, we have that
$\exp(\phi_2-\frac{\eta_2}{\eta_1}\phi_1)\in\ot_\xsa(U_{j'})$, which
contradicts \eqref{eq_psi21}. Hence
$\exp(\phi_2)\notin\ot_\xsa(U_{j'})$.
\proofend



\begin{cor}\label{cor_twisted_ram_exp_case}
  Let $l\in\integergz$ $\omega,\phi_1,\phi_2\in\zcz$, such that
  $-v(\omega)>\underset{j=1,2}\max\{\frac{-v(\phi_j)}l+1\}$. Let $S$
  an open sector of amplitude $2\pi$, $\zeta$ an inverse of $z\mapsto
  z^l$ def\mbox{}ined on $S$. The following conditions are equivalent.  
\ben
\item $\phi_1\circ\zeta\neq\phi_2\circ\zeta$ .
\item At least one of the following two conditions is verif\mbox{}ied:
\ben
\item there exists $U\in\opxsac$, $U\subset S$ such that
  $\exp(\phi_{1}\circ\zeta+\omega)\in\ot_\xsa(U)$ and
  $\exp(\phi_{2}\circ\zeta+\omega)\notin\ot_\xsa(U)$;
\item there exists $V\in\opxsac$, $V\subset S$ such that
  $\exp(\phi_{1}\circ\zeta+\omega)\notin\ot_\xsa(V)$ and
  $\exp(\phi_{2}\circ\zeta+\omega)\in\ot_\xsa(V)$.  
\een 
\een

\end{cor}

\emph{Proof.}
\emph{(ii)$\Rightarrow$(i)}. Obvious.

\emph{(i)$\Rightarrow$(ii)}. Set $\mu_l(z):=z^l$. 

Suppose now that
$\phi_1\circ\zeta\neq\phi_2\circ\zeta$. It follows that, for any
$\lambda\in\reagz$, $\lambda(\phi_1+\omega\circ\mu_l)\neq\phi_2+\omega\circ\mu_l$. Consider the sector $\zeta(S)$ of amplitude
$\frac{2\pi}l$. Since $-v(\omega)\geq2$, then
$\frac{2\pi}{l}>-\frac{2\pi}{lv(\omega)}$. Hence by Proposition \ref{prop_exp_case}
 there exists $\zeta(U)\subset\zeta(S)$ such
that either $\exp(\phi_1+\omega\circ\mu_l)\in\ot_\xsa(\zeta(U))$ and
$\exp(\phi_2+\omega\circ\mu_l)\notin\ot_\xsa(\zeta(U))$ or
viceversa. By Proposition \ref{prop_pullback_ot}, it follows that
either $\exp(\phi_1\circ\zeta+\omega)\in\ot_\xsa(U)$ and
$\exp(\phi_2\circ\zeta+\omega)\notin\ot_\xsa(U)$ or viceversa.  

\proofend

\begin{rem}
{\rm There is another way to prove Proposition
  \ref{prop_exp_case}. Let us brief\mbox{}ly summarize it.

Let $\phi_1,\phi_2\in\zcz\setminus\{0\}$, such that, for any
$\lambda\in\reagz$, $\phi_1\neq\lambda \phi_2$. Let $n_k=-v(\phi_k)$ ($k=1,2$).

Let $C_{\phi_k,A}$ be the boundary of the set $U_{\phi_k,A}$ ($k=1,2$). One has
 that $C_{\phi_k,A}$ is the set of the zeros of a polynomial
 $Q_{\phi_k,A}(x,y)\in\rea[x,y]$ ($k=1,2$). Further $C_{\phi_k,A}$ has $2n_k$ distinct
 branches at $0$ determined by the Puiseux's series
 $\sigma_{\phi_k,A,1}(x),\ldots,\sigma_{\phi_k,A,2n_k}(x)$ obtained by solving the equation
 $Q_{\phi_k,A}(x,y)=0$ ($k=1,2$) with respect to $y$.

One checks that the f\mbox{}irst $n_k$ terms of $\sigma_{\phi_k,A,j}(x)$
do not depend on $A$ ($k=1,2$, $j=1,\ldots,2n_k$). Further, it turns out that there
exists $\theta_k\in[0,2\pi[$ such that the tangent at $0$ of the graph
of $\sigma_{\phi_k,A,j}(x)$ has slope
$\tan\bigr{\theta_k+j\frac\pi{2n_k}}$ ($k=1,2$, $j=1,\ldots,2n_k$).

If $n_1\neq n_2$ or $\theta_1\neq\theta_2$, the result follows easily.

If $n_1=n_2$ and $\theta_1=\theta_2$, one checks that there exist
$\overline j\in\bigc{1,\ldots,2n_1}$ and $r\in\bigc{1,\ldots,n_1}$ such
that the $r$-th coef\mbox{}f\mbox{}icients of
$\sigma_{\phi_1,A,\overline j}(x)$ and $\sigma_{\phi_2,A,\overline
  j}(x)$ are dif\mbox{}ferent. Hence, there are inf\mbox{}initely many
relatively compact subanalytic open sets concentrated along some
$\theta_1+j\frac\pi{2n_1}$ f\mbox{}itting between
$\sigma_{\phi_1,A,\overline j}(x)$ and $\sigma_{\phi_2,A,\overline
  j}(x)$. Choosing $U$ among these sets, one obtains that one
exponential is tempered on $U$ and the other is not.  

This procedure is more intuitive than the proof we chose to expose
here but it is more technical and much longer.  
}
\end{rem}

\section{Tempered solutions and $\Omega$-f\mbox{}iltered local
  systems} \label{SEC_FORM_INV}

In the f\mbox{}irst part of this section we are going to prove that the
tempered solutions induce a fully faithful functor on good models. In
the second part we will prove that the datum of tempered solutions of
a meromorphic connection $\mc M$ is equivalent to the data of
determinant polynomials and holomorphic solutions of $\mc M$. 

Let us recall that M. Kashiwara, in \cite{rims}, proves that,
given a complex analytic manifold $X$ and an object $\mc M$ of the
bounded derived category of $\D_X$-modules with regular holonomic
cohomology, 
$$  R\Ho[\rho_!\D_X]{\rho_!\M}{\otxsa}
\simeq 
R\Ho[\rho_!\D_X]{\rho_!\M}{\mathcal O_X} \ . $$

Given a complex analytic curve $Y$, $x_0\in Y$, $\M\in\mod_h(\D_{*x_0})
$, there exists a neighborhood $X\subset Y$ of $x_0$ such that $\mc M$
is a holonomic $\D_X$-module. By chosing a local coordinate $z$ near
$x_0$, we can suppose that $X=\com$ and $x_0=0$. Recall that for
$\omega\in\zcz$, we set $\mc L^\omega:=\D_\com\exp(\omega)$. 

Set 
\begin{eqnarray*}
 \sh\M  &  :=  &  \Ho[\D_X]{\M}{\mathcal O_X} \ , \\
 \st\M  &  :=  &  \Ho[\rho_!\D_X]{\rho_!\M}{\otxsa} \ ,\\
 \ms S^t_\omega(\M)  &  :=  &  \Ho[\rho_!\D_X]{\rho_!(\M\otimes\mc L^{\omega})}{\otxsa} \ ,
\end{eqnarray*}
the functors $\rho_*$ and $\rho_!$ being def\mbox{}ined in Subsection
\ref{SUBSEC_RECALL_XSA}. 

\subsection{Tempered solutions and non-ramif\mbox{}ied $\Omega$-graduations}

The main results of this subsection are Proposition
\ref{prop_tensor_case} and Theorem \ref{thm_ff} below. First,
let us describe explicitely the subanalytic sheaf $\ms S^t\bigr{\mc
  L^\phi\otimes\mc R}$, for $\phi\in\zcz$ and $\mc R$ a regular
holonomic $\D_{*0}$-module. 



Recall the def\mbox{}inition of the sets $U_{\phi,A}$ given in \eqref{eq_UphiA}.

\begin{lem}\label{lem_sheaves_st_sh}
Let $\phi\in\zcz$, $\mc R$ a regular holonomic $\D_{*0}$-module. Then
$$\st{\mc L^\phi\otimes \mc R}\simeq\underset{A>0}{\varinjlim}\ \rho_*\ \sh{\mc R}_{U_{\phi,A}}$$ 
\end{lem}

\emph{Proof.} If $\phi=0$, the result follows from the fact that $\ms
S^t(\mc R)\simeq\ms S(\mc R)$. 

Suppose $\phi\neq0$. Let $V\in\opxsac$ be connected and simply connected. If $0\in V$, then
clearly $\Gamma\bigr{V,\st{\mc L^\phi\otimes\mc
    R}}\simeq0$. Otherwise, the $\com$-vector space
$\Gamma\bigr{V,\sh{\mc L^\phi\otimes\mc R}}$ 
has f\mbox{}inite dimension $r$ and is generated by
$h_1(z)\exp(\phi(z)),$ $\ldots,$ $h_r(z)\exp(\phi(z))$, for 
$h_1,\ldots,h_r\in\mathcal O_\com(V)$, such that there exist $C,M>0$ satisfying
$$ \phantom{a}\qquad\qquad\qquad C|z|^M\leq|h_j(z)|\leq\bigr{C|z|^M}^{-1} \qquad (z\in V,\,j=1,\ldots,r) \ .$$

In particular, since $\Gamma\bigr{V,\st{\mc L^\phi\otimes\mc
    R}}\simeq\Gamma\bigr{V,\sh{\mc L^\phi\otimes\mc R}}\cap\otxsa(V)$
we have that  
$$ \Gamma\bigr{V,\st{\mc L^\phi\otimes\mc R}}\simeq\left\{
\begin{array}{ll}
\Gamma\bigr{V,\sh{\mc L^\phi\otimes\mc R}} &  \text{ if }
\exp(\phi)\in\otxsa(V) \\
0 & \text{ otherwise}
\end{array}
\right .  \ . $$



The conclusion follows by Proposition \ref{prop_generalization_arx}.

\proofend

We can now state and proof

\begin{prop}\label{prop_tensor_case}
Let $\phi_1,\phi_2\in\zcz$, $\phi_2\neq0$, $\mc R_1,\mc R_2$ regular holonomic
$\D_{*0}$-modules. 
If, for any $\lambda\in\reagz$, $\phi_1\neq\lambda\phi_2$ 
$$  \hob[\com_{\xsa}]{\st{\mc L^{\phi_1}\otimes \mc R_1}}{\st{\mc
    L^{\phi_2}\otimes \mc R_2}}
\simeq0 \ .$$
Otherwise,
\begin{equation}
  \label{eq_functorially}
  \hob[\com_{\xsa}]{\st{\mc L^{\phi_1}\otimes \mc R_1}}{\st{\mc
    L^{\phi_2}\otimes \mc R_2}}
\simeq\hob[\D_{*0}]{\mc R_1}{\mc R_2}
\end{equation}
functorially in $\mc R_1,\mc R_2$.
\end{prop}

\emph{Proof.} Suppose that, for any $\lambda\in\reagz$,
$\phi_1\neq\lambda\phi_2$. By Proposition \ref{prop_exp_case}\emph{(ii)}, there
exists $W\in\opxsac$ such that
\begin{enumerate}
\item there exists $A_0>0$ such that, for any $A\geq A_0$, $W\subset U_{\phi_1,A}$,
\item for any $B>0$, $W\not\subset U_{\phi_2,B}$.
\end{enumerate}

In particular, for any $A>A_0$ and $B>0$,
\begin{equation}\label{eq_A>B}
U_{\phi_1,A}\not\subset U_{\phi_2,B} \ .
\end{equation} 
Combining \eqref{eq_A>B} and the fact that $\sh{\mc R_1}$ and $\sh{\mc R_2}$ are locally
constant sheaves on $\com^\times$, we obtain, for any $A>A_0$ and $B>0$,
$$ \hob[\com_X]{\sh{\mc R_1}_{U_{\phi_1,A}}}{\sh{\mc R_2}_{U_{\phi_2,B}}}=0 \ . $$
Now, using Lemma \ref{lem_sheaves_st_sh}, we obtain
\begin{gather*}
\hob[\com_\xsa]{\st{\mc L^{\phi_1}\otimes \mc R_1}}{\st{\mc L^{\phi_2}\otimes \mc R_2}} 
\simeq \hspace{80mm}\phantom{a}\\
  \begin{align*}
    \hspace{40mm}\simeq & \ 
    \underset{A>0}{\varprojlim}\underset{B>0}{\varinjlim}\hob[\com_X]{\sh{\mc R_1}_{U_{\phi_1,A}}}{\sh{\mc R_2}_{U_{\phi_2,B}}} \\ 
= & \  0 \ .
  \end{align*}
\end{gather*}


Suppose now there exists
$\lambda\in\reagz$ such that $\phi_1=\lambda\phi_2$. Then
\begin{equation}\label{eq_Uphi1=Uphi2}
U_{\phi_1,\lambda A}= U_{\phi_2,A} \ .
\end{equation} 

We need the following

\begin{lem}\label{lem_UphiA_hom_C}
Let $\phi\in\zcz\setminus\{0\}$ and $A\in\reagz$. The set
 $$U_{\phi,A}:=\{z\in\com^\times;\ \re\phi(z)<A\} $$
is homotopically equivalent to $\com^\times$.
\end{lem}

\emph{Proof.} We prove the result in three steps: $\phi=1/z$,
$\phi=1/z^n$ and $\phi\in\zcz$.

First suppose that $\phi(z)=\frac1z$. Then $U_{\phi,A}$ is the complementary of a closed
ball and the result is obvious. 

Suppose now that $\phi(z)=\frac1{z^n}$, for some $n\in\integergz$. Let
$\mu_n:\com\to\com$, $z\mapsto z^n$. Then
$U_{\phi,A}=\mu_n^{-1}(U_{\frac1z,A})$ and the conclusion follows.  

Suppose now that $\phi\in\zcz$ and $-v(\phi)=n$. Mimicking the proof
of Proposition \ref{prop_+-temp}, there exists a  biholomorphism
$\eta$ such that $\phi(\eta(z))=\frac1{z^n}$. The conclusion follows. 
\nopagebreak
\proofend

Let us conclude the proof of Proposition \ref{prop_tensor_case}.

We have the following sequence of isomorphisms
\begin{gather*}
\hob[\com_\xsa]{\st{\mc L^{\phi_1}\otimes \mc R_1}}{\st{\mc
    L^{\phi_2}\otimes \mc R_2}} \simeq \qquad\qquad\qquad\qquad\qquad\qquad\phantom{a}
   \\  
\begin{align*}
\phantom{a}\qquad\qquad\qquad\qquad\qquad&\simeq  \underset{A>0}{\varprojlim}\,\underset{B>0}{\varinjlim}\,\hob[\com_X]{\sh{\mc R_1}_{U_{\phi_1,A}}}{\sh{\mc R_2}_{U_{\phi_2,B}}} \\
& \simeq   \underset{A>0}{\varprojlim}\,\underset{B>0}{\varinjlim}\,\hob[\com_X]{\sh{\mc R_1}_{U_{\phi_1,A}}}{\sh{\mc R_2}_{U_{\phi_1,A}}} \\
& \simeq   \underset{A>0}{\varprojlim}\,\underset{B>0}{\varinjlim}\,\hob[\com_X]{\sh{\mc R_1}_{\com^\times}}{\sh{\mc R_2}_{\com^\times}} \\
& \simeq   \hob[\D_{*0}]{\mc R_1}{\mc R_2}\ ,
\end{align*}
\end{gather*}

where the f\mbox{}irst isomorphism follows from Lemma
\ref{lem_sheaves_st_sh}, the second from \eqref{eq_Uphi1=Uphi2} and the third
from Lemma \ref{lem_UphiA_hom_C}.

The conclusion follows.

\proofend




We can now state the main results of this subsection.

\begin{thm}
Let $\ou{\phi\in\Sigma_1}\oplus\mc L^{\phi}\otimes\mc R_\phi$ and 
$\ou{\psi\in\Sigma_2}\oplus\mc L^{\psi}\otimes\mc P_\psi$ be two good
models. The following conditions are equivalent. 
\begin{enumerate}
\item 
$ \ms S^t\Bigr{\ou{\phi\in\Sigma_1}\oplus\mc L^{\phi}\otimes\mc R_\phi}
_{X\setminus\{0\}}
\simeq
\ms S^t\Bigr{\ou{\psi\in\Sigma_2}\oplus\mc L^{\psi}\otimes\mc P_\psi}
_{X\setminus\{0\}} $.
\item 
There exist $\phi_1,\ldots,\phi_d\in\zcz$ such that,
$\ou[d]{j=1}{\coprod}\reagz\phi_j=\reagz\Sigma_1=\reagz\Sigma_2$ and,
for any $j=1,\ldots,d$, $\ou{\phi\in\Sigma_1\cap\reagz\phi_j}{\oplus}\mc
R_\phi\simeq\ou{\phi\in\Sigma_2\cap\reagz\phi_j}{\oplus}\mc P_\psi$. 
\end{enumerate} 
\end{thm}

\emph{Proof.}
$(ii)\Rightarrow(i)$ Combining Proposition \ref{prop_tensor_case} and
the fact that $\ms S^t(\cdot)_{X\setminus\{0\}}$ is fully faithfull on
the category of regular holonomic $\D_{*0}$-modules, we have that  
\begin{gather*}
  \mathrm {Hom}_{\com_\xsa}
\Bigr{
\ms S^t\Bigr{\ou{\phi\in\Sigma_1}\oplus\mc L^{\phi}\otimes\mc R_\phi}_{X\setminus\{0\}},
\ms S^t\Bigr{\ou{\psi\in\Sigma_2}\oplus\mc L^{\psi}\otimes\mc P_\psi}_{X\setminus\{0\}}
}
\simeq\qquad\qquad\qquad\\
\simeq\ou{\phi\in\Sigma_1}{\oplus}\ou{\psi\in\Sigma_2}{\oplus}\hob[\D_{*0}]{\ms
  S^t(\mc L^\phi\otimes \mc R_\phi)}{\ms S^t(\mc L^\psi\otimes\mc
  P_\psi)} \\
\qquad\qquad\qquad\simeq\ou[d]{j=1}{\oplus}\ \ou{\phi\in\Sigma_1\cap\reagz\phi_j}{\oplus}\hob[\D_{*0}]{\ms
  S^t(\mc L^\phi\otimes \mc R_\phi)}{\ou{\psi\in\Sigma_2\cap\reagz\phi_j}{\oplus}\ms S^t(\mc L^\psi\otimes\mc P_\psi)} \\
\simeq\ou[d]{j=1}{\oplus}\hob[\D_{*0}]{\ou{\phi\in\Sigma_1\cap\reagz\phi_j}{\oplus}\mc
  R_\phi}{\ou{\psi\in\Sigma_2\cap\reagz\phi_j}{\oplus}\mc P_\psi} \ . 
\end{gather*}
The functoriality of \eqref{eq_functorially} allows to conclude.

$(i)\Rightarrow(ii)$ First let us suppose that
$\reagz\Sigma_1\neq\reagz\Sigma_2$. Hence either
$\reagz\Sigma_1\not\subset\reagz\Sigma_2$ or
$\reagz\Sigma_2\not\subset\reagz\Sigma_1$. Suppose the latter. 

There exists $\psi\in\Sigma_2$ such that for any $\phi\in\Sigma_1,\lambda\in\reagz$, $\psi\neq\lambda\phi$.

Suppose that $\psi\neq0$. By Proposition \ref{prop_tensor_case}, we have

$$  \hob
[\com_{\xsa}]
{\st{\ou{\phi\in\Sigma_1}{\oplus}\mc L^{\phi}\otimes \mc R_\phi}_{X\setminus\{0\}}}
{\st{\mc L^{\psi}\otimes \mc P_\psi}_{X\setminus\{0\}}}
\simeq0 \ .$$

It follows that

$$\st{\ou{\phi\in\Sigma_1}{\oplus}\mc L^{\phi}\otimes \mc R_\phi}_{X\setminus\{0\}}
\not\simeq
\st{\ou{\psi\in\Sigma_2}{\oplus}\mc L^{\psi}\otimes \mc P_\psi}_{X\setminus\{0\}}
 \ . $$

Suppose that $\psi=0$, then $0\notin\Sigma_1$ and, by Proposition \ref{prop_tensor_case},

$$  \hob
[\com_{\xsa}]
{\st{\mc L^{\psi}\otimes \mc P_\psi}_{X\setminus\{0\}}}
{\st{\ou{\phi\in\Sigma_1}{\oplus}\mc L^{\phi}\otimes \mc R_\phi}_{X\setminus\{0\}}}
\simeq0 \ .$$

It follows that

$$\st{\ou{\phi\in\Sigma_1}{\oplus}\mc L^{\phi}\otimes \mc R_\phi}_{X\setminus\{0\}}
\not\simeq
\st{\ou{\psi\in\Sigma_2}{\oplus}\mc L^{\psi}\otimes \mc P_\psi}_{X\setminus\{0\}} \ .
$$

The case $\reagz\Sigma_1\not\subset\reagz\Sigma_2$ is treated similarly.

Now let us suppose that
$\reagz\Sigma_1=\reagz\Sigma_2=\ou[d]{j=1}{\coprod}\reagz\phi_j$ and,
there exists $j'\in\{1,\ldots,d\}$, such that  
$$
\ou{\phi\in\Sigma_1\cap\reagz\phi_{j'}}{\oplus}\mc R_\phi
\not\simeq
\ou{\psi\in\Sigma_2\cap\reagz\phi_{j'}}{\oplus}\mc P_\psi \ .
$$

By Proposition \ref{prop_tensor_case} we have that
\begin{gather*}
  \hob[\com_\xsa]{\ms S^t\Bigr{\ou{\phi\in\Sigma_1}\oplus\mc
      L^{\phi}\otimes\mc R_\phi}_{X\setminus\{0\}}}{\ms
    S^t\Bigr{\ou{\psi\in\Sigma_2}\oplus\mc L^{\psi}\otimes\mc
      P_\psi}_{X\setminus\{0\}}} 
\simeq\qquad\qquad\qquad\\
\qquad\qquad\qquad\qquad\qquad\qquad\qquad\simeq\ou[d]{j=1}{\oplus}\hob[\D_{*0}]{\ou{\phi\in\Sigma_1\cap\reagz\phi_j}{\oplus}\mc
  R_\phi}{\ou{\psi\in\Sigma_2\cap\reagz\phi_j}{\oplus}\mc P_\psi} \ . 
\end{gather*}
The functoriality of \eqref{eq_functorially} allows to conclude.

\proofend

We conclude the study of tempered solutions of good models  with

\begin{thm}\label{thm_ff}
 Let $\omega\in\zcz$, $-v(\omega)\geq k$. The functor 
$$\begin{array}{rrcl}
  \ms S^t_\omega(\cdot):  &  \msf{GM}_k  &  \lra  &  \mod(\com_{\xsa})\\
& \M & \longmapsto & \Ho[\rho_!\D_X]{\rho_!(\M\otimes\mc L^{\omega})}{\otxsa} \ ,$$ 
\end{array}$$
is fully faithful.
\end{thm}

\emph{Proof.} Clearly it is suf\mbox{}f\mbox{}icient to prove that, given
$\phi_1,\phi_2\in\zcz$, $k>\max\{-v(\phi_1),-v(\phi_2)\}$, $\mc
R_1,\mc R_2$ regular holonomic $\D_{*0}$-modules, the functor of
tempered solutions induces the isomorphism 
\begin{equation}\label{eq_ff}  
\hob[\D_{*0}]{\mc L^{\phi_1}\otimes \mc R_1}{\mc
    L^{\phi_2}\otimes \mc R_2}
\simeq 
\hob[\com_{\xsa}]{\ms S^t_\omega\bigr{\mc L^{\phi_1}\otimes \mc R_1}}{\ms S^t_\omega\bigr{\mc
    L^{\phi_2}\otimes \mc R_2}}
 \ .
\end{equation}

Let us prove \eqref{eq_ff}. 

First, suppose that $\phi_1\neq\phi_2$. Then 
$$\hob[\D_{*0}]{\mc L^{\phi_1}\otimes\mc R_1}{\mc L^{\phi_2}\otimes\mc
  R_2}=0  \ . $$ 
Moreover, as, for any $\lambda\in\reagz$,
$\lambda(\phi_1+\omega)\neq\phi_2+\omega$, Proposition
\ref{prop_tensor_case} implies that,
$$ \hob[\com_{\xsa}]{\ms S^t_\omega\bigr{\mc L^{\phi_1}\otimes \mc R_1}}{\ms S^t_\omega\bigr{\mc
    L^{\phi_2}\otimes \mc R_2}}=0 \ .
$$
 
Now, suppose that $\phi_1=\phi_2$. The result follows from Proposition
\ref{prop_tensor_case} and the fact that 
$$  \hob[\D_{*0}]{\mc R_1}{\mc R_2}\simeq\hob[\D_{*0}]{\mc
  L^{\phi_1}\otimes \mc R_1}{\mc  L^{\phi_2}\otimes \mc R_2}  \ . $$ 

\proofend

\subsection{Tempered solutions of ordinary dif\mbox{}ferential equations}

We begin this subsection by proving the analogue of Lemma
\ref{lem_sheaves_st_sh} in the case of ramif\mbox{}ied determinant
polynomials i.e. with non-integer exponents.  

Recall that, for $Y\subset X$, $Y_{\xsa}$ is the subanalytic
site on $Y$ induced by $\xsa$, in particular the open sets of $Y_{\xsa}$
are of the form $U\cap Y$ for $U\in\opxsac$. For $\mc F\in\mod(k_\xsa)$,
we denote by $\mc F|_Y$ the restriction of $\mc F$ to $Y_\xsa$.

Recall the def\mbox{}initions of $\Omega(\mc M)$ and $r_{\phi,\mc M}$
(resp. $U_{\phi,\ep}$) given in Def\mbox{}inition \ref{df_formal_invariants}
(resp. Corollary \ref{cor_generalization_arx_ram}).

\begin{lem}\label{lem_st_conn}
Let $\mc M\in\mod_h(\D_{*0})$, $\theta\in\rea$,
$Y:=X\setminus(\rea_{\geq0} e^{i\theta})$. 
Then
$$  \st {\mc M}\big|_Y\simeq \ou{\phi\in\Omega(\mc M)}{\oplus}\ \underset{\ep>0}{\varinjlim}\
\rho_*\ \com^{r_{\phi,\mc M}}_{Y,U_{\phi,\ep}}  \ . $$
\end{lem}

\emph{Proof.} As $\mc M$ is f\mbox{}ixed, for sake of simplicity, we drop the
index $\mc M$ in the symbol $r_{\phi,\mc M}$.

Let $V\in\Op(\yxsa)$ connected. By the Hukuhara-Turrittin's Asymptotic
Theorem \ref{thm_asymptotic_HT}, the $\com$-vector space $\sh {\mc
  M}(V)\subset\mathcal O(V)$ is generated by $\bigc{h_{\phi,j}\exp(\phi)}_{
\begin{subarray}{l}
\phi\in\Omega(\mc M)\\ j\in\{1,\ldots,r_\phi\}
\end{subarray}
}$.

Hence
\begin{eqnarray*}
\st {\mc M}(V)  & \simeq  &  \sh{\mc M}(V)\cap\ot(V) \\
  &\simeq  &
  \Bigc{\sum_{\phi\in\Omega(\mc M)}\sum_{j=1}^{r_\phi}c_{\phi,j}h_{\phi,j}\exp(\phi)\in\ot(V);\
    c_{\phi,j}\in\com} \ .
\end{eqnarray*}
Since, for $\phi\in\Omega(\mc M), j\in\{1,\ldots,r_\phi\}$,
$h_{\phi,j}\exp(\phi)$ are $\com$-linearly independent functions and
$h_{\phi,j},h_{\phi,j}^{-1}\in\ot(V)$, one has that 
$$ \sum_{\phi\in\Omega(\mc M)}\sum_{j=1}^{r_\phi}c_{\phi,j}h_{\phi,j}\exp(\phi)\in\ot(V) $$
if and only if $\exp(\phi)\in\ot(V)$ for $c_{\phi,j}\neq0$.

The conclusion follows.
\proofend

We denote by $\msf{for}$ the functor from the category of
$\Omega$-f\mbox{}iltered local systems to the category of local
systems on $S^1$. 

\begin{thm}\label{thm_final}
Let $k\in\integer_{>0}$, $\mc M_1,\mc M_2\in\mod_h(\D_{*0})_k$ and 
$\omega\in\zcz$ such that $-v(\omega)>k$. The following conditions
are equivalent. 
\begin{enumerate}
\item $\ms S^t_\omega\bigr{\mc M_1}_{X\setminus\{0\}}\simeq\ms S^t_\omega\bigr{\mc M_2}_{X\setminus\{0\}}$. 
\item 
  \begin{enumerate} 
    \item $\msf{for}\bigr{\ms S^\Omega(\mc
      M_1)}\simeq\msf{for}\bigr{\ms S^\Omega(\mc M_2)}$ and 
    \item for any $\theta\in S^1$, $\ms S^\Omega\bigr{\mc
     M_1}_\theta\simeq\ms S^\Omega\bigr{\mc M_2}_\theta$ as
  $\Omega_\theta$-graded $\com$-vector spaces. 
 \end{enumerate} 
\end{enumerate}
\end{thm}

\emph{Proof.}
\emph{(i)$\Rightarrow$(ii)}. First remark that $\msf{for}\bigr{\ms
  S^\Omega(\mc M_1)}\simeq\msf{for}\bigr{\ms S^\Omega(\mc M_2)}$ if
and only if $\sh{\mc M_1}_{X\setminus\{0\}}\simeq\sh{\mc
  M_2}_{X\setminus\{0\}}$. Since $\rho^{-1}\st{\mc
  M_j}_{X\setminus\{0\}}\simeq\sh{\mc M_j}_{X\setminus\{0\}}$, the
condition \emph{(a)} is proved. 

Suppose now that there exists $\theta\in S^1$ such that $\ms
S^\Omega(\mc M_1)_\theta\not\simeq\ms S^\Omega(\mc M_2)_\theta$. Then,
either $\Omega(\mc M_1)\neq\Omega(\mc M_2)$ or there exists
$\phi\in\Omega(\mc M_1)\cap\Omega(\mc M_2)$ such that $r_{\phi,\mc
  M_1}\neq r_{\phi,\mc M_2}$. In the former case, combining the ideas
of the f\mbox{}irst part of the proof of Proposition \ref{prop_tensor_case}
with Lemma \ref{lem_st_conn} and Corollaries
\ref{cor_generalization_arx_ram} and \ref{cor_twisted_ram_exp_case},
we obtain that, for any $\theta\in\rea$, $\ms S^t_\omega(\mc
M_1)|_{X\setminus\rea_{\geq0}e^{i\theta}}\not\simeq\ms
S^t_\omega(\mc M_2)|_{X\setminus\rea_{\geq0}e^{i\theta}}$. In the
latter case the result follows easily from Lemma \ref{lem_st_conn}. 

\emph{(ii)$\Rightarrow$(i)}. 
Set $\ms S_\omega(\cdot):=\ms S(\cdot\otimes\mc L^\omega)$ ($j=1,2$). 

Let $\theta_1,\theta_2\in\rea$, $\theta_1\neq\theta_2\
(\mathrm{mod}\,2\pi)$, $Y_j:=X\setminus\rea_{\geq0} e^{i\theta_j}$
($j=1,2$).

Since for any $\theta\in S^1$, $\ms
S^\Omega\bigr{\mc M_1}_\theta\simeq\ms S^\Omega\bigr{\mc M_2}_\theta$,
then $\Omega(\mc M_1)=\Omega(\mc M_2)$ and $r_{\phi,\mc
  M_1}=r_{\phi,\mc M_2}$. In particular, Lemma \ref{lem_st_conn}
implies that   
$$ \ms S^t_\omega(\mc M_1)|_{Y_1}\simeq  \ms S^t_\omega(\mc M_2)|_{Y_1}$$
$$ \ms S^t_\omega(\mc M_1)|_{Y_2}\simeq  \ms S^t_\omega(\mc M_2)|_{Y_2} \ .$$

Now, since 
$\msf{for}\bigr{\ms S^\Omega(\mc M_1)}
\simeq
\msf{for}\bigr{\ms S^\Omega(\mc M_2)}$, we have that
$ \sh{\mc M_1}_{X\setminus\{0\}}
  \simeq
  \sh{\mc M_2}_{X\setminus\{0\}}  $
which implies $\ms S_\omega\bigr{\mc M_1}_{X\setminus\{0\}}\simeq\ms S_\omega\bigr{\mc M_2}_{X\setminus\{0\}}$.
We conclude thanks to the commutative
diagram below. Roughly speaking, it says that $\ms S^t_\omega(\mc
M_j)_{X\setminus\{0\}}$ is completely determined by $\ms
S^t_\omega(\mc M_j)|_{Y_1}$, $\ms S^t_\omega(\mc M_j)|_{Y_2}$ and $\ms S_\omega(\mc M_j)_{X\setminus\{0\}}$
($j=1,2$). 

$$ 
\xymatrix{
0\ar[r]& \ms S^t_\omega(\mc
M_j)_{X\setminus\{0\}}\ar[r]\ar@{^{(}->}[d] &  \ms S^t_\omega(\mc
M_j)_{Y_1}\oplus\ms S^t_\omega(\mc M_j)_{Y_2}\ar[r] \ar@{^{(}->}[d] &
\ms S^t_\omega(\mc M_j)_{Y_1\cap Y_2}\ar@{^{(}->}[d]\ar[r] & 0\\ 
0\ar[r]& \ms S_\omega(\mc M_j)_{X\setminus\{0\}}\ar[r] &  \ms
S_\omega(\mc M_j)_{Y_1}\oplus\ms S_\omega(\mc M_j)_{Y_2}\ar[r] &  \ms
S_\omega(\mc M_j)_{Y_1\cap Y_2}\ar[r] & 0 
} 
$$
\nopagebreak[4]\proofend


\vspace{10mm}

{\small{\sc Dipartimento di Matematica Pura ed Applicata,

Universit{\`a} di Padova,

Via Trieste 63, 35121 Padova, Italy.}

E-mail address: $\texttt{gmorando@math.unipd.it}$}

\end{document}